\documentclass[11pt,a4paper]{amsart}

\usepackage{amsmath}
\usepackage{amssymb}
\usepackage{latexsym}
\usepackage{xypic}
\usepackage{amscd,amsthm}

\newcommand{\dual}{\makebox[0mm]{}^{{\scriptstyle\vee}}}

\usepackage{graphicx}
\usepackage{psfrag}
\usepackage{epic}
\usepackage{eepic}
\input{epsf}

\newcommand{\card}[1]{{\mid\! #1 \!\mid}}
\newcommand{\pro}[2]{\langle #1, #2 \rangle}

\newcommand{\qu}[1]{{\overline{#1}}}

\def\lra{{\,\Leftrightarrow\,}}
\def\ra{{\,\Rightarrow\,}}

\def\dt{\overline{d}}
\def\Pt{{\widetilde{P}}}

\def\cydim{{\rm CY\-\!\!\!-\!\!\-dim}}
\def\cydimgross{{\rm CY\-\!\!-\!\-dim}}

\def\hst{h_\text{st}}
\def\Est{{E_{\rm st}}}
\def\St{{\widetilde{S}}}

\def\p{{\tilde{p}}}

\def\Dt{{\tilde{\Delta}}}
\def\C{{\mathbb C}}
\def\N{{\mathbb N}}

\def\P{{\mathbb P}}
\def\Q{{\mathbb Q}}
\def\R{{\mathbb R}}
\def\Z{{\mathbb Z}}

\def\dual{{\vee}}
\def\msd{{m_{\sigma^\dual}}}
\def\ns{{n_\sigma}}
\def\MRq{{\overline{M}_\R}}
\def\MRt{{\widetilde{M}_\R}}
\def\Mt{{\widetilde{M}}}
\def\Nt{{\widetilde{N}}}
\def\Mq{{\overline{M}}}
\def\NRq{{\overline{N}_\R}}

\def\NRt{{\Nt_\R}}

\def\Nq{{\overline{N}}}
\def\Nqs{{\Nq^{(r)}}}

\def\nt{{\tilde{\nabla}}}

\def\F{{\mathcal F}}

\def\calP{{\mathcal P}}

\def\T{{\mathcal T}}

\def\V{{\mathcal V}}

\def\rk{{\rm rk}}

\def\conv{{\rm Conv}}
\def\aff{{\rm aff}}
\def\lin{{\rm lin}}
\def\pos{{\R_{\geq 0}}}

\def\Vol{{\rm Vol}}

\def\rank{{\rm rank}}

\def\deg{{\rm deg}}

\def\Hom{{\textup{Hom}}}

\def\NR{{N_{\R}}}
\def\NQ{{N_{\Q}}}
\def\MR{{M_{\R}}}
\def\MQ{{M_{\Q}}}
\def\Nq{{\qu{N}}}
\def\Mq{{\qu{M}}}

\def\rand{\partial}

\def\intr{{\rm int}}

\newtheorem*{theorem*}{Theorem}
\newtheorem*{corollary*}{Corollary}
\newtheorem{theorem}{Theorem}[section]
\newtheorem{definition}[theorem]{Definition}
\newtheorem{remark}[theorem]{Remark}
\newtheorem{example}[theorem]{Example}
\newtheorem{corollary}[theorem]{Corollary}
\newtheorem{proposition}[theorem]{Proposition}
\newtheorem{lemma}[theorem]{Lemma}
\newtheorem{conjecture}[theorem]{Conjecture}
\newtheorem*{acknowledgment}{Acknowledgment}
\newtheorem{question}[theorem]{Question}

\begin{document}
\title{Combinatorial aspects of mirror symmetry}

\author[Victor Batyrev]{Victor Batyrev}
\address{Mathematisches Institut, Universit\"at T\"ubingen, 
Auf der Morgenstelle 10, 72076 T\"ubingen, Germany}
\email{ victor.batyrev@uni-tuebingen.de}

\author[Benjamin Nill]{Benjamin Nill}
\address{Arbeitsgruppe Gitterpolytope, Freie Universit\"at Berlin, 
Arnimallee 3, 14195 Berlin, Germany}
\email{ nill@math.fu-berlin.de}

\begin{abstract}
The purpose of this paper is to review some combinatorial ideas behind 
the   mirror symmetry for Calabi-Yau 
hypersurfaces and complete intersections in Gorenstein toric Fano varieties. 
We suggest as a basic combinatorial object the notion of a Gorenstein polytope 
of index $r$. A natural combinatorial duality for $d$-dimensional 
Gorenstein polytopes of index $r$ extends the well-known polar duality 
for reflexive polytopes (case $r=1$).
We consider the Borisov duality 
between two  nef-partitions as a duality between  
two Gorenstein polytopes $P$ and $P^*$ of index $r$ together with  
selected   special $(r-1)$-dimensional simplices $S \subset P$ and 
$S' \subset P^*$. Different choices of these simplices suggest 
an interesting relation  to  Homological Mirror Symmetry.  

\vspace{-0.5ex}
\end{abstract}

\maketitle

\section*{Introduction}

Several papers of the first author and  Borisov  
 were devoted to a  combinatorial method for  
explicit constructions of mirror pairs of Calabi-Yau varieties 
(see \cite{Bor93, Bat94, BB96a, BB97}). 
 A starting point of this method was the notion of a 
{\em reflexive polytope} introduced in \cite{Bat94}. 
A reflexive polytope $P \subset \R^d$ is a convex $d$-dimensional 
polytope with vertices in $\Z^n$ containing the origin $0 \in \R^d$ 
in its interior such that its dual (polar) polytope 
\[ P^* := \{ y =(y_1, \ldots, y_d)  \in \R^d\; : \; 
\sum_{k=1}^d x_k y_k \geq -1 \;\; \forall x =(x_1, \ldots, x_d) \in P \}. 
\]
is again a lattice polytope (i.e., all vertices of $P^*$  
belong to $\Z^n$). In this case $P^*$ is also a reflexive 
polytope and $(P^*)^* = P$. For any convex lattice polytope 
$\Delta \subset \R^d$ we denote by $X_{\Delta}$  a generic hypersurface 
in $(\C^*)^d$ defined by the equation 
\[ F_{\Delta}(z) := \sum_{m \in \Z^d \cap \Delta} a_m z^m =0, \] 
where the coefficients $a_m$ $(m \in \Z^d \cap P)$ are independent 
variables and $z=(z_1, \ldots, z_d)$ are standard complex coordinates 
on  $(\C^*)^d$. If $P$ is reflexive, 
then  $X_P \subset(\C^*)^d $ admits a natural 
compactification $\overline{X}_P$ which is 
a $(d-1)$-dimensional Calabi-Yau variety. Moreover, the duality between 
$P$ and $P^*$ becomes the mirror duality between $(d-1)$-dimensional 
Calabi-Yau varieties $\overline{X}_P$ and
$\overline{X}_{P^*}$

Many explicit constructions of mirrors of Calabi-Yau complete intersections
suggested by the first author and van Straten 
in \cite{BS95} motivated  Borisov to generalize the polar 
duality 
for reflexive polytopes to a more general duality for 
so called  
{\em nef-partitions} (see  \cite{Bor93}). From a combinatorial point of view, 
a {\em nef-partition of length $r$} is a decomposition 
$P = P_1 + \cdots + P_r$ 
of a $d$-dimensional reflexive polytope 
$P \subset \R^n$ into a Minkowski sum of $r$ lattice polytopes 
$P_1, \ldots, P_r$ such that $0 \in P_i$ for all $i =1, \ldots, r$ 
(it is not assumed that each $P_i$ has the maximal dimension $d$).  
It can be shown that for any $j=1, \ldots, r$ the polytope
\[ Q_j =  \{ y   \in \R^d\; : \; 
\sum_{k=1}^d x_k y_k \geq -\delta_{ij} \;\; \forall x 
\in P_i , \;\forall i =1, \ldots, r \}
\]
has vertices in $\Z^d$, $0 \in Q_j$ and the Minkowski sum 
$Q = Q_1 + \cdots + Q_r$ is again a $d$-dimensional reflexive polytope. 
Moreover, the lattice polytopes $P_1, \ldots, P_r$ 
can be obtained from lattice polytopes $Q_1, \ldots, Q_r$ 
by the same formula: 
\[ P_i =  \{ y   \in \R^d\; : \; 
\sum_{k=1}^d x_k y_k \geq -\delta_{ji} \;\; \forall x 
\in Q_j, \;\forall j =1, \ldots, r \}.
\]
The Minkowski sum decomposition $Q = Q_1 + \cdots + Q_r$ is called 
the {\em dual nef-partition}. Some examples of nef-partitions and their 
duals are given in Examples \ref{nefp}, \ref{dual-nefp} and 
\ref{nef-p2}.  If 
 $X_{\{P_i\}}$ is a generic complete intersection 
in $(\C^*)^d$ defined by generic equations
\[ F_{P_1}(z) = \cdots = F_{P_r}(z)  =0, \]   
then   $X_{\{P_i\}}$  admits a natural compactification 
$\overline{X}_{\{P_i\}}$ which is 
a $(d-r)$-dimensional Calabi-Yau variety.  Moreover, the duality of  
nef-partitions 
$\{P_i \}$ and $\{Q_j\}$ becomes the mirror duality between 
$(d-r)$-dimensional Calabi-Yau varieties  
$\overline{X}_{\{P_i\}}$ and
$\overline{X}_{\{Q_j\}}$. 

If $P \subset \R^n$ is a reflexive polytope, $P = P_1 + \cdots + 
P_r$ is a nef-partition and $Q = Q_1 + \ldots + 
Q_r$ is the dual nef-partition, then the lattice polytopes $P_1, \ldots, 
P_r$ satisfy some additional combinatorial conditions: 

\begin{itemize}
\item[1.] $\pos{(P_k)} \cap \pos{(P_l)} = 0$ for all $1 \leq k < l \leq r$;

\item[2.] the convex hull ${\conv}\, \{ P_1, \ldots, P_r \}$ 
is the dual reflexive polytope $Q^*$;

\item[3.] $Q^* \cap \Z^d = (P_1 \cap \Z^d) \cup \cdots \cup 
(P_r \cap \Z^d)$. Moreover, $m \in Q^* \cap \Z^d$ is a vertex 
of $Q^*$ if and only if $m$ is a nonzero vertex of $P_i$ for 
some $i \in \{1, \ldots, r\}$.
  
\end{itemize}

For lattice polytopes 
$\Delta_1, \ldots, \Delta_r \in \R^d$ we denote by $\Delta_1 * 
\cdots * \Delta_r$ the polytope in  $\R^d \times \R^r$ 
which is the convex hull 
\[ {\rm Conv}\{ (\Delta_1 \times e_1), \ldots, (\Delta_r \times e_r) 
\}, \]
where $e_1, \ldots, e_r$ is the standard basis of $\R^r$. We call 
$\Delta_1 * 
\cdots * \Delta_r$  {\em Cayley polytope associated with} 
$\Delta_1, \ldots, \Delta_r$.  

Let  $P$ and $Q$ be two $d$-dimensional reflexive polytopes and 
$P = P_1 + \cdots +P_r$ and $Q = Q_1 + \cdots + Q_r$ nef-partitions 
which are  dual to each 
other. We consider two $(d+r-1)$-dimensional polytopes  
 $\tilde{P} := P_1 * \cdots * P_r$, 
$\tilde{Q} := Q_1 * \cdots * Q_r$ and two  $(d+r)$-dimensional cones  
\[ C_{\tilde{P}}: = \R_{\geq 0}(P_1 \times e_1)  + 
\cdots + \R_{\geq 0}(P_r \times e_r) \subset \R^d \times \R^r, \]
\[ C_{\tilde{Q}} : = \R_{\geq 0}(Q_1 \times e_1)  + 
\cdots +  \R_{\geq 0}(Q_r \times e_r) \subset \R^d \times \R^r. \]
It was proved in \cite[Theorem 4.6]{BB97} that the  
cones  $C_{\tilde{P}}$ and  $C_{\tilde{Q}}$ are dual to each other 
with respect to standard scalar product on $\R^d \times \R^r$. 
Moreover, there exists a 
unique  lattice  point $p \in  r \tilde{P}$ (resp.  $q \in  r \tilde{Q}$) 
in the relative interior of $r \tilde{P}$ (resp. $r \tilde{Q}$) 
such that $ r \tilde{P} - p$ (resp. $ r \tilde{Q} - q$) is a 
$(d+r-1)$-dimensional reflexive polytope. 
For this reason, the duality for nef-partitions of length $r$ 
can be considered as a special case of the more general duality for 
so called   {\em reflexive Gorenstein cones of index $r$}. 
The duality for reflexive Gorenstein cones gives new  possibilities 
for constructing  mirrors of Calabi-Yau varieties. In particular, it 
allows to construct mirrors of rigid-Calabi-Yau $3$-folds 
in form of so called {\em generalized Calabi-Yau varieties} or 
{\em Landau-Ginzburg orbifolds} \cite{CDP93,BB97}. 

Recall that a  $(d+r-1)$-dimensional lattice polytope $P \subset 
\R^{d+r-1}$ is called a {\em Gorenstein polytope of index $r$}, 
if $rP$ contains a single interior lattice point $p$ and 
$rP - p$ is a reflexive polytope. In this case, the $(d+r)$-dimensional 
cone 
\[ C_P:= \R_{\geq 0} (P \times 1) \subset \R^{d+r-1} \times \R. \]
is called {\em reflexive Gorenstein cone of index $r$} 
associated with $P$. Let  
$${C}_P^\dual := \{ y \in \R^d \times \R, \; \sum_k {x_k}{y_k} \geq 0, 
\;\forall x \in C_P \}$$ be the dual cone.
Then there exists another Gorenstein polytope $P^*$ 
of index $r$ such that $C_{P^*}$ isomorphic (via the action  of 
$GL(d+r, \Z)$) to the cone ${C}_P^\dual$. The polytope $P^*$ is 
uniquely determined up to the action of  the affine group 
$AGL(d+r-1, \Z)$ on  $\R^{d+r-1}$. It is called 
the {\em dual Gorenstein polytope of index $r$}.   
In this paper, we want to look at the duality of  reflexive Gorenstein 
cones of index $r$ from the viewpoint of the duality of  
Gorenstein polytopes of index $r$. For $r=1$ the latter coincides with 
the already known polar duality for reflexive polytopes.  
It is important to consider an  additional structure on Gorenstein 
polytopes $P$ which is defined by a choice of a special $(r-1)$-dimensional
lattice simplex $S \subset P$. The combinatorial notion of a 
{\em special simplex} was recently 
introduced in \cite{Ath05} in order to prove a conjecture of Stanley. 
If a Gorenstein polytope $P$ of index 
$r$ and its dual $P^*$ both contain special $(r-1)$-simplices, we natually 
obtain two  nef-partitions dual to each other. Moreover, the existence 
of special simplices in $P$ and $P^*$ can be considered as 
another characterization of nef-partitions and their duals (see 
Proposition \ref{splittheo}). The choice of special 
$(r-1)$-simplices is not unique. Therefore, one may obtain many different
nef-partitions from the same 
pair  of Gorenstein polytopes $(P,P^*)$. We expect that 
different choices of special 
$(r-1)$-simplices $S \subset P$ and $S' \subset P^*$
define  Calabi-Yau complete
intersections which are equivalent 
from the viewpoint of Homological Mirror Symmetry.  
\medskip

The paper is organized as follows. In Section 1 we recall the notions 
of Gorenstein polytopes, reflexive Gorenstein cones and their dualities.
In Section 2 we establish a bijective correspondence  
between Cayley structures of length $r$ on a Gorenstein polytope $P$ 
of index $r$ and special $(r-1)$-simplices in the dual Gorenstein 
polytope $P^*$. In Section 3 we review the notion of nef-partitions
and their duality as well as combinatorial properties of nef-partitions. 
In Section 4 we define 
the $E_{\rm st}$-function  of a Gorenstein polytope $P$ and formulate some 
open questions and conjectures. 
Finally, in Section 5 we give a precise formulation of a conjecture  
which connects  special $(r-1)$-simplices in a Gorenstein 
polytope $P$ with the Homological Mirror Symmetry. 
In the last Section 6 we consider 
as an addendum to Section 3 some combinatorial operations with 
nef-partitions.

\begin{acknowledgment}{\rm 
We are grateful to anonymous referees for valuable remarks. 
The second author thanks Christian Haase for discussions, motivating questions and 
suggestions.
}
\end{acknowledgment}

\section{Gorenstein polytopes and reflexive Gorenstein cones}

\subsection{Gorenstein polytopes}

Let $M \cong \Z^d$ be a lattice of rank $d$, $N = \Hom_\Z(M,\Z)$ 
the dual lattice, $\pro{\cdot}{\cdot}\,:\, M \times N \to \Z$ the natural 
pairing. 
We define $\MQ := M \otimes_\Z \Q$, $\MR := M \otimes_\Z \R$, 
$\NQ := N \otimes_\Z \Q$, $\NR:= N \otimes_\Z \R $. 
By a polytope $\Delta \subset \MR$ we always mean a 
 convex hull of finitely many points in $\MR$. 
The set of vertices of $\Delta$ is  denoted by $\V(\Delta)$, 
the relative interior by $\intr(\Delta)$, and  
the boundary by $\rand \Delta$. For arbitrary subset  
$A \subseteq \MR$ we denote by $\conv(A)$ (resp. $\aff(A)$)  
the convex (resp. affine) hull of $A$. 

From now on we assume $\Delta \subseteq \MR$ to be a $d$-dimensional 
{\em lattice polytope}, i.e., $\V(\Delta) \subset M$.

\smallskip

\begin{definition}{\rm Assume that  $\Delta$ has the origin $0 \in M$ 
in its interior. Then 
the {\em dual polytope} of $\Delta$ is defined as 
\[\Delta^* := \{y \in \NR \;:\; \pro{y}{x} \geq -1 \;\forall\, x \in \Delta\}.\]
This is a rational polytope (i.e., $\V(\Delta^*) \subseteq \NQ$)  
containing the origin $0 \in N$ in its interior. 

If $\Delta$ contains an interior lattice point $m \in M$ 
(not necessary $0$), we say that 
$(\Delta-m)^*$ is the {\em dual polytope of $\Delta$ with respect to $m$}.
}\end{definition}

\begin{definition}{\rm Let $\Delta \subseteq \MR$ be a 
lattice polytope containing the origin $0 \in M$ in its interior. 
We say $\Delta$ is {\em reflexive}, if the dual polytope 
$\Delta^*$ is also a lattice polytope. 
In this case, $0 \in \Delta$ is a single  interior lattice point.

In general, if a lattice polytope $\Delta \subset \MR$ 
contains a single interior lattice 
point $m$ (not necessary $0$), we say that 
$\Delta$ is {\em reflexive with respect to $m$}, if $\Delta-m$ is reflexive.
}
\end{definition}

\begin{remark}{\rm It is easy to see that 
a lattice polytope $\Delta \subseteq \MR$ is 
reflexive if and only if 
every facet of $\Delta$ has integral distance $1$ 
from its unique interior lattice point. 
Recall that a point $x \in \MR$ is said to have 
{\em integral distance} $\delta$ from a facet $F$ of $\Delta$, if 
there is a primitive lattice point $u \in N$ 
such that $u \perp F$ and $\pro{u}{F} - \pro{u}{x} = \delta$.
}
\end{remark}

\begin{remark}{\rm 
The following example shows that in dimension $ \geq 3$ there exist 
lattice polytopes containing $0$ as a single interior lattice point which 
 are not reflexive. Let $\Delta$ be the $3$-dimensional lattice polytope 
with vertices $\pm (1,0,0)$, $\pm (0,1,0)$, $\pm (1,1,2)$. 
Then $(0,0,0) \in \Delta$ is a single  interior lattice point, however, 
the integral distance from $0$ to the facet  
$$F := \conv((1,0,0),(0,1,0),(1,1,2))$$ is $2$, since 
for the primitive vector $u := (2,2,-1)^t$ we have  $u \perp F$ and 
$\pro{u}{F} = 2$.
}
\end{remark}

\begin{definition}{\rm Let $r$ be a positive integer. 
A lattice $d$-dimensional polytope $\Delta \subset \MR$ is called 
{\em Gorenstein of index} $r$, if $r \Delta$ contains an interior 
lattice point $m$ and $r\Delta -m$ is reflexive. 
}
\end{definition}

\begin{example} 
{\rm Let $\Delta \subset \R^d$ be the standard unit cube: $ 0 \leq x_i \leq 1$
$\forall i = 1, \ldots, d$. Then $2\Delta - (1, \ldots, 1)$ is reflexive. 
Therefore $\Delta$ is a Gorenstein polytope of index $2$.
}
\end{example}

It is clear that Gorenstein polytopes can be characterized as follows:

\begin{remark}{\rm 
A lattice polytope $\Delta$ is a Gorenstein polytope 
of index $r$ if and only if 
there is a  rational 
point $x \in \intr(\Delta) \cap \frac{1}{r}M$ having the integral distance 
$1/r$ from any facet of $\Delta$.
}
\end{remark}

\smallskip

\subsection{Reflexive Gorenstein cones}

Let $\Mq$ and $\Nq$ be lattices of rank $\dt$ which 
are dual to each other. 
Let $\sigma \subset \MRq$ be a $\dt$-dimensional rational finite 
polyhedral 
cone with vertex $0 \in \Mq$. Then 
 \[\sigma^\dual := \{y \in \NRq \,:\; \pro{x}{y} \geq 0 \;
\forall\, x \in \sigma\}\]
is called the 
{\em dual cone} of $\sigma$.  
It is also a $\dt$-dimensional rational finite polyhedral cone with 
vertex $0$.

We recall some definitions from  \cite{BB97}:

\begin{definition}{\rm  A $\dt$-dimensional rational finite 
polyhedral 
cone $\sigma$ is called {\em Gorenstein cone}, 
if it is generated by finitely many lattice points 
which are  contained in an affine 
hyperplane $\{x \in \MRq \,:\, \pro{x}{n} = 1\}$ for some $n \in \Nq$. 
In this case, the lattice point $n \in \intr(\sigma^\dual)$ is 
uniquely determined, and we denote it  by 
$\ns$. One has $ \intr(\sigma^\dual) \cap \Nq = \ns + \sigma^\dual \cap \Nq$. 
We define the $k$th slice of $\sigma$ as 
$\sigma_{(k)} := \sigma \cap \{x \in \MRq \,:\, \pro{x}{\ns} = k\}$. 

}
\label{Gor-cone}
\end{definition}

\begin{definition}{\rm 
Let $\sigma \subseteq \MRq$ be a Gorenstein cone. Then the lattice 
polytope $\Dt := \sigma_{(1)}$ 
is called the {\em support} of $\sigma$; it is a  $(\dt-1)$-dimensional 
lattice polytope with respect to the 
 affine lattice $\Mt := \Mq \cap \{x \in \MRq \,:\, \pro{x}{\ns} = 1\}$. 

Let  $\Delta \subseteq \MR$ be  a lattice polytope of dimension $d$. 
We set  $\Mq := M \oplus \Z$, $\dt := d+1$,  
$\Mt := M \times \{1\}$, $\Dt := \Delta \times \{1\}$. 
Then $\pos(\Dt)$ is a Gorenstein cone with support $\Dt$.
}
\end{definition}

Using the two  constructions above, 
we obtain  a $1$-to-$1$ correspondence between 
$\dt$-dimensional Gorenstein cones $\sigma \subset \Mq$ 
considered up to action of $GL(\dt, \Z)$ on 
$\Mq$ and $(\dt-1)$-dimensional  lattice polytopes  $\Dt \subset \Mt_\R$
considered up 
to an affine transformation from $AGL(\dt-1, \Z)$.

\begin{definition}{\rm 
A Gorenstein cone $\sigma$ is called {\em reflexive}, 
if $\sigma^\dual$ is also a Gorenstein cone. 
In this case, we denote by $\msd$ the unique lattice point in $\Mq$ 
such  that the affine hyperplane $\{y \in \NRq \,:\, \pro{\msd}{y} = 1\}$ 
contains all lattice points generating the cone  
$\sigma^\dual$.  The number 
\[r = \pro{\msd}{\ns}\]
is called  the {\em index}  of the reflexive Gorenstein cone 
$\sigma$. 
Using   $(\sigma^\dual)^\dual = \sigma$, 
we see that $r$ is also the index of the reflexive Gorenstein cone 
$\sigma^\dual$. Thus, we obtain a natural duality of 
 reflexive Gorenstein cones of index $r$.
}
\label{refl-G}
\end{definition}

It was observed in \cite[Prop.2.11]{BB97} that  
reflexive Gorenstein cones correspond to Gorenstein polytopes:

\begin{proposition}
Let $\sigma$ be a $\dt$-dimensional Gorenstein cone with support $\Dt$. 
Then the following statements are equivalent:

\begin{enumerate}
\item $\sigma$ is a reflexive Gorenstein cone of index $r$;
\item $\sigma_{(r)}$ is a reflexive polytope;
\item $\Dt$ is a Gorenstein polytope of index $r$
\end{enumerate}
In this case, $\sigma_{(r)}$ $($resp. $(\sigma^\dual)_{(r)}$$)$ is a 
reflexive polytope with the unique interior lattice point 
$\msd$ (resp. $\ns$). 
\label{reflcone}
\end{proposition}

The  duality of reflexive Gorenstein cones of index $r$ defines 
a natural duality of Gorenstein polytopes of index $r$:

\begin{definition}{\rm Let $\Delta \subseteq \MR$ be a $d$-dimensional 
Gorenstein polytope of index $r$. 
Then the {\em dual Gorenstein polytope} $\Delta^*$ is defined as the 
support of the dual cone 
$\sigma^\dual$, where $\sigma$ is the reflexive Gorenstein 
cone associated to $\Delta$. By \ref{reflcone}, 
$\Delta^*$ is also a $d$-dimensional Gorenstein 
polytope of index $r$.
}
\end{definition}

\begin{remark} 
{\rm We want to compare the duality of $d$-dimensional Gorenstein polytopes 
of index $r$ with the polar duality of reflexive polytopes. 
If $r =1$, then the dualities are the same. For any $r \geq 1$ 
two Gorenstein polytopes  $\Delta$ and $\Delta^*$ are combinatorially 
dual to each other. 
However, if $r >1$, then the two reflexive polytopes $r\Delta$ and $r\Delta^*$  are 
{\em not} dual to each other as reflexive polytopes. 
In order to 
obtain the dual reflexive polytope  $(r\Delta)^*$ from the Gorenstein polytope 
$\Delta^*$ one must replace the lattice $N$ by another one  (see 
\ref{refine} and \ref{refldual} below). } 
\end{remark}

\smallskip
Let us consider the duality of Gorenstein polytopes of index $r$ 
in more detail.

\smallskip \subsection{Duality of Gorenstein polytopes}
\label{detail}

\begin{definition}{\rm 
Let $\sigma \subseteq \MRq$ be a Gorenstein cone 
of index $r$ with support $\Dt \subseteq \MRt$.
We define two lattices 
\[\Mq_{(r)} := \{x \in \Mq \,:\, \pro{x}{\ns} = 0 
\;(\text{mod}\, r)\}, \text{ and } \Nq^{(r)} := \Nq + \frac{1}{r}\Z \ns.\]
The lattice $\Mq_{(r)}$ is dual to $\Nq^{(r)}$ and one has
 $\Mq_{(r)} \subseteq \Mq$ and $\Nq \subseteq \Nq^{(r)}$. 
\label{refine}}
\end{definition}

The cone $\sigma$ considered with respect to the new lattice 
$\Mq_{(r)}$ is a Gorenstein cone of index $1$, and the reflexive polytope
$\sigma_{(r)} = r \Dt$ is its support. Thus,  we obtain: 

\begin{proposition} Let $\Dt \subseteq \MRt$ be a Gorenstein polytope of index 
$r$ and  $\Dt^* \subset  \NRt$ be the dual Gorenstein polytope. 
If  $\sigma \subseteq \MRq$ is a reflexive Gorenstein cone 
of index $r$ with support $\Dt \subseteq \MRt$,  
then the dual reflexive polytope $(r \Dt)^*$ equals 
$\Dt^*$ with respect to the refined affine lattice 
$$\Nq^{(r)} \cap \{y \in \NR \,:\, \pro{\msd}{y} = 1\}.$$
\label{refldual}
\end{proposition}

\begin{proposition}
Let $\rho$ be a facet of $\sigma^\dual$. We denote by $\lin(\rho)$ 
the minimal linear subspace in $\NRq$ containing  $\rho$.
Then
\[\lin(\rho) \cap \Nq = \lin(\rho) \cap \Nqs.\]
In particular, one has 
\[\rand \sigma^\dual \cap \Nq = \rand \sigma^\dual \cap \Nqs \]
and 
\[   \rand \Dt^* \cap \Nt =  \rand (r \Dt)^* \cap \Nq^{(r)}. \]
\label{boundarylemma}
\end{proposition}

\begin{proof}

Let $x \in \lin(\rho) \cap \Nqs$. Since $\rho$ is 
a facet of $\sigma^\dual$, there is a vertex $v$ of $\Dt$ 
such that $\pro{v}{y} = 0$ for all $y \in \rho$. 
From $\pro{v}{\ns} = 1$ we get 
$\Nq = \Z \ns \oplus (\Nq \cap v^\perp)$. 
Hence $\Nqs = \frac{1}{r}\Z \ns \oplus (\Nq \cap v^\perp)$. Now,  
$\pro{v}{x} = 0$ implies $x \in \Nq$. 
\end{proof}

\smallskip

\smallskip

\begin{remark} 
{\rm From the combinatorial point of view,  
the duality of Gorenstein polytopes ``coincides'' 
with the polar duality. However, 
there exists a  ``subtlety'' coming  from the occuring lattices. 
Let $\Dt$ be a Gorenstein polytope of index  $r > 1$.  Then 
the  Gorenstein polytope $\Dt^*$ has no interior lattice points. 
Proposition  \ref{boundarylemma} shows that the reflexive polytope 
$(r \Dt)^*$ and $\Dt^*$ have the same sets of boundary lattice points. 
Assume that $r$ is not a prime number. 
Let $k > 1$ be any proper divisor of $r$. Then $k \Dt$ is a 
Gorenstein polytope of index $r/k$. 
The dual Gorenstein polytope $(k \Dt)^*$ has the same lattice points 
on the boundary as $(r \Dt)^*$ and   $\Dt^*$. Therefore, 
two Gorenstein polytopes  $(k \Dt)^*$ and  $\Dt^*$ have the same 
combinatorial structure and also the same lattice point structure. However, 
$(k \Dt)^*$ and $\Dt^*$ are Gorenstein polytopes of different index, 
so they are {\em not} isomorphic as lattice polytopes.} 
\end{remark}

Consider an example of this phenomenon:

\begin{example}{\rm
 Let $\Delta$ be a  unimodular $3$-dimensional simplex, i.e., it has 
lattice volume $1$. Then  $\Delta$ 
is a Gorenstein polytope of index $4$. We have $\Delta \cong \Delta^*$, i.e. 
$\Delta$ is a selfdual Gorenstein polytope. The dual reflexive polytope 
$(4 \Delta)^*$ is a  $3$-dimensional reflexive simplex with $4$ lattice points 
on the boundary. Now consider  $2\Delta$ as a Gorenstein polytope of index $2$.
The dual Gorenstein polytope $(2 \Delta)^*$ 
is a Gorenstein simplex of index $2$ 
 having lattice volume $2$ and having only $4$ vertices 
as lattice points  (the same number of lattice points as 
in $\Delta^*$). 

The polytope $4 (\Delta^*)$ is a reflexive simplex with 
$34$ lattice points on the boundary. The polytope 
$2 (2 \Delta)^*$ is a self-dual reflexive simplex with 
$10$ lattice points on the boundary. 
Pictures  of these polytopes are given  in Example 
\ref{fig1} where   
the polytope $2 \Delta$ is denoted by $\Dt^*$.
}
\end{example}

\section{Cayley polytopes and special simplices}

In this section we consider  Cayley polytopes of length $r$ 
and give a criterion for a Gorenstein polytope $\Dt$ of index $r$ 
 to be a Cayley polytope of length  $r$. For this purpose we use 
the notion of a special $(r-1)$-dimensional simplex in the dual 
Gorenstein polytope $\Dt^*$.

\subsection{Cayley polytopes} 

Let   $M,N$ be dual lattices of rank $d$; 
and let $\Delta_1, \ldots, \Delta_r \subseteq \MR$ be lattice polytopes.

\begin{definition}{\rm
Consider  the lattice $\Mq := M \oplus \Z^r$, where $e_1, \ldots, e_r$ is 
the standard lattice basis of $\Z^r$. Then $\Mq$ contains the affine 
sublattice $\Mt := M \oplus (\aff(e_1, \ldots, e_r) \cap \Z^r)$.
We identify  $\Delta_i \subseteq \MR$ with  
$\Dt_i := \Delta_i \times e_i \subseteq \MRt$ for $i = 1, \ldots, r$.
The polytope  $\Delta_1 * \cdots * \Delta_r := 
\conv(\Dt_1, \ldots, \Dt_r) \subseteq \MRt$ 
is called the {\em Cayley polytope of length $r$} 
associated to $\Delta_1, \ldots, \Delta_r$. The cone 
$$\pos(\Delta_1 * \cdots * \Delta_r) = \pos \Dt_1 + \cdots + \pos \Dt_r  
\subseteq \Mq$$ is called the associated {\em Cayley cone}.
\label{standard}
}
\end{definition}

If $\aff(\Delta_1, \ldots, \Delta_r) = \MR$, 
then the associated Cayley cone is a Gorenstein cone of dimension $d+r$, and 
$\Dt = \Delta_1 * \cdots * \Delta_r$ 
is its support.

\smallskip

\begin{definition}{\rm 
Let $e_1, \ldots, e_r \in \Z^r$ be the standard  lattice basis. 
Then $S_{r-1} := \conv(e_1, \ldots, e_r)$ is called 
an {\em unimodular simplex} of dimension $r-1$.
}
\end{definition}

\begin{proposition}
Let $\sigma \subseteq \MRq$ be a Gorenstein cone with support 
$\Dt \subseteq \MRt$. Then the following 
statements are equivalent:

\begin{enumerate}
\item $\sigma$ is a Cayley cone associated to $r$ lattice polytopes;
\item $\Dt$ is a Cayley polytope of length $r$;
\item There is a lattice projection $\Mt \twoheadrightarrow \Z^r$, 
which maps $\Dt$ surjectively on  $S_{r-1}$;
\item There are nonzero $e^*_1, \ldots, e^*_r \in \sigma^\dual \cap \Nq$, 
such that 
$e^*_1 + \cdots + e^*_r = \ns.$
\end{enumerate}
Moreover, the lattice vectors  $e^*_1, \ldots, e^*_r$ in 
$(4)$ form a  part of a basis of $\Nq$ and the Cayley structure of $\Dt$ 
is uniquely determined by $r$ polytopes 
\[\Dt_i := \{x \in \Dt \;:\; \pro{x}{e^*_j} = 0 
\text{ for } j \not= i \} \; (i=1, \ldots, r).\]
These polytopes  have properties 
$\pro{\Dt_i}{e^*_i} = 1$, and $\Dt = \conv(\Dt_1, \ldots, \Dt_r)$.
\label{char}
\end{proposition}

\begin{proof}

(1), (2), (3) are obviously equivalent.

(1) $\ra$ (4): Let   $b_1, \ldots, b_d, e_1, \ldots, e_r$ be a basis of 
the lattice $\Mq$, where $b_1, \ldots, b_d$ is a basis of the  lattice 
$M$ and  
$e_1, \ldots, e_r$ the standard basis of $\Z^r$. 
Consider the  dual basis $b^*_1, \ldots, b^*_d, e^*_1, \ldots, e^*_r$ of 
the lattice $\Nq$.
Then we have $\pro{x_j}{e^*_i} = \delta_{i j}$ 
for $i,j \in \{1, \ldots, r\}$ and $x_j \in \Dt_j$ ($\delta_{i j}$ 
denotes the Kronecker symbol, i.e., $\delta_{i j} = 1$, for $i = j$, 
and $0$ otherwise). 
Therefore we obtain  $e^*_1, \ldots, e^*_r \in \sigma^\dual$, 
$\pro{\Dt}{e^*_1 + \cdots + e^*_r} = 1$. The latter implies 
$e^*_1 + \cdots + e^*_r = \ns$.

(4) $\ra$ (1): Define the polytopes 
$\Dt_i := \{x \in \Dt \;:\; \pro{x}{e^*_j} = 0$ for 
$j \not= i$$\}$ for $i=1, \ldots, r$. 
Let $x$ be a vertex of $\Dt$. Then $1 = \pro{x}{e^*_1 + \cdots + e^*_r}$. 
Since 
$0 \leq \pro{x}{e^*_i} \in \Z$ for $i = 1, \ldots, r$, we get that 
$x \in \Dt_i$ for some $i$ and $\pro{x}{e_i^*} =1$. Hence $\Dt = 
\conv(\Dt_1, \ldots, \Dt_r)$. 
Take any $i \in \{1, \ldots, r\}$. We claim that $\Dt_i \not= \emptyset$. 
Indeed, otherwise we  would have $\pro{\Dt}{e_i^*} =0$. 
This would contradict the condition 
$e_i^* \neq 0$.   
Hence the 
mapping $x \to (\pro{x}{e^*_1}, \ldots, \pro{x}{e^*_r})$ gives 
a projection of lattices $\Mt \twoheadrightarrow \Z^r$ 
such that  the image of  $\Dt$ is  $S_{r-1}$. In particular,  
 $e^*_1, \ldots, e^*_r$ form a part of a basis of  $\Nq$.

\end{proof}

Recall the following definition  from \cite[Def. 3.9]{BB97}:

\begin{definition}{\rm 
A reflexive Gorenstein cone $\sigma$ of index $r$ is said to be 
{\em completely split}, 
if $\sigma$ is the Cayley cone associated to $r$ lattice polytopes.
}
\end{definition}

By \ref{char}, 
we obtain the following  characterization of completely split 
reflexive Gorenstein cones:  

\begin{corollary}
Let $\sigma \subseteq \MRq$ be a reflexive Gorenstein cone of index $r$ 
with support $\Dt$. Then the following 
statements are equivalent:
\begin{enumerate}
\item $\sigma$ is completely split;
\item there exist lattice points  $e^*_1, \ldots, e^*_r \in \Dt^* 
\cap \Nq$ such that  $$e^*_1 + \cdots + e^*_r = \ns.$$
\end{enumerate}
\smallskip
Moreover, there is a $1$-to-$1$ correspondence between all possible 
collections 
$\{ e^*_1, \ldots, e^*_r \}$  of $r$ lattice points in $\Dt^*$ 
such that  $e^*_1 + \cdots + e^*_r =\ns$ and all possible 
Cayley polytope structures of length $r$ of the Gorenstein polytope $\Dt$.
\label{compl}
\end{corollary}

\begin{proof}
(1) $\ra$ (2): Assume that $\sigma$ is completely split.  By \ref{char}(4) 
there exist lattice points  $e^*_1, \ldots, e^*_r \in \sigma^\dual$ 
such that  $e^*_1 + \cdots + e^*_r = \ns$. Since $r = \pro{\msd}{\ns} = 
\pro{\msd}{e^*_1} + \cdots + \pro{\msd}{e^*_r}$ and 
 $1 \leq \pro{\msd}{e^*_i} \in \Z$ for all $i = 1, \ldots, r$, we get 
$\pro{\msd}{e^*_i} = 1$ for all $i = 1, \ldots, r$. Thus, all 
lattice points  $e^*_1, \ldots, e^*_r$ belong to  $\Dt^*$. 

(2) $\ra$ (1): The statement follows immediately from \ref{char}, because  
 the index of $\sigma$ equals $\pro{\msd}{\ns} = 
\pro{\msd}{e^*_1} + \cdots + \pro{\msd}{e^*_r} = r$.
\end{proof}

By \ref{char}, every completely split reflexive Gorenstein cone 
$\sigma \subset \MRq$ determines 
$r$ lattice polytopes $\Dt_1, \ldots, \Dt_r$.   
On the other hand, one can start with  $r$ lattice polytopes $\Delta_1, \ldots, 
\Delta_r \subset \MR$ and ask whether the corresponding 
Cayley cone is reflexive of index $r$, or whether 
the Cayley polytope $\Delta_1 * \cdots * \Delta_r$  is a Gorenstein polytope 
of index $r$. The answer  is contained in  \cite[Prop.3.6]{BB97}:

\begin{theorem}
Let  $\Delta_1, \ldots, \Delta_r \subset \MR$ be lattice polytopes such that 
the Minkowski sum $\Delta_1 + \cdots + \Delta_r$  has dimension $d$.
Then the  following statements are equivalent:
\begin{enumerate}
\item The associated Cayley cone $\sigma$ is reflexive of index $r$;
\item $\Delta_1 * \cdots * \Delta_r$ is a Gorenstein polytope of index $r$;
\item $\Delta_1 + \cdots + \Delta_r$ is a reflexive polytope 
(with interior lattice point $m$).
\end{enumerate}
In this case, the lattice point  $\msd  \in \sigma$ (see Def. \ref{refl-G}) 
equals $m \times (e_1 + \cdots + e_r)$.
\label{sum}
\end{theorem}

For convenience, we give  a purely  combinatorial proof of this statement. 

\begin{proof} (1) $\lra$ (2): Follows from Prop. \ref{reflcone}.

(2) $\ra$ (3): Let $\Dt= \Delta_1 * \cdots * \Delta_r$ be a Gorenstein 
polytope of index $r$. Denote by $\pi$ the projections $\Mq \to \Z^r$ and 
$\MRq \to \R^r$. Since $r\Dt$ is reflexive, there exists a unique  
lattice point $m \in \intr (r\Dt)$. Moreover, the $\pi$-image of 
$m$ can only be the unique interior point $s_1:= (1,\ldots, 1)$ in the 
reflexive simplex $r S_{r-1}$ (where $S_r = \conv(e_1, \ldots, e_r)$). We 
remark that the intersection of the affine subspace $\pi^{-1}(s_1)$ with 
$r\Dt$ is exactly $\Delta_1 + \ldots + \Delta_r$. Since $m$ has integral
distance $1$ from any facet of $r\Dt$, the lattice point $m$  has also 
integral distance $1$ from 
any facet of  $\Delta_1 + \ldots + \Delta_r$ (facets of  
$\Delta_1 + \ldots + \Delta_r$ are intersections of facets of  $r\Dt$ with 
  the affine subspace $\pi^{-1}(s_1)$).

(3) $\ra$ (1): 
Assume that the lattice polytope $\Delta:= \Delta_1 + \cdots + 
\Delta_r$ is reflexive. Denote by $m$ the unique interior 
lattice point in $P$. Consider the Cayley cone $\sigma$ 
over $\Dt= \Delta_1 * \cdots * \Delta_r$ and denote
by $\msd$ the lattice point $m \times (e_1 + \cdots + e_r)$ in 
$\intr(\sigma)$. We have to show that $\pro{\msd}{z} =1$ for 
any primitive lattice vector $z \in \Nq$ generating a $1$-dimensional
face $\pos z$ of the dual cone $\sigma^\dual \subset \NRq$. 

Let $\Gamma := z^{\perp} \cap \sigma$ the facet of $\sigma$ 
which is orthogonal to $z$. There are the following two possibilities:

Case 1. There exists $i \in \{ 1, \ldots, r \}$ such that $\Gamma \cap 
(\Delta_i 
\times e_i) = \emptyset$. Then $\Gamma$ is orthogonal to $e_i^*$, because
$\Gamma$ is generated by some lattice points from the polytopes 
$(\Delta_j \times e_j)$ $(j =1, \ldots, r)$ and $\pro{x}{e_i^*} = \delta_{ij}$ 
for $x \in (\Delta_j \times e_j)$. Since $e_i^*$ is primitive and $e_i^* 
\in \sigma^\dual$, we get 
$z=e_i^*$ and $\pro{\msd}{z} = \pro{e_1 + \cdots + e_r}{e_i^*} =1$. 

Case 2. For any  $i \in \{ 1, \ldots, r \}$ there exists a lattice point
$p_i \in \Delta_i \times e_i$ which is contained in $\Gamma$. 
In this case, 
$\Gamma$ is a Cayley cone of $r$ lattice polytopes     
$P_i:= \Gamma \cap (\Delta_i \times e_i)$ $ (i =1, \ldots, r)$ such that 
$p_i \in P_i$ $( 1 \leq i \leq r)$.   
Moreover, 
$P_{\Gamma}:= \Gamma \cap (\Delta \times (e_1 + \cdots + e_r)) = 
P_1 + \cdots + P_r$ is a $(d-1)$-dimensional 
face of $\Delta \times (e_1 + \cdots + e_r) \cong \Delta_1 + \cdots + \Delta_r$.

Let $M'$ be the $(d+1)$-dimensional sublattice in $\Mq$ defined 
by $r-1$ equations $\pro{x}{e_1^*} = \cdots = \pro{x}{e_r^*}$. 
Then $\sigma_\Delta:= \MR' \cap \sigma$ is exactly 
the cone over the $d$-dimensional reflexive  
polytope $\Delta = \Delta_1 + \cdots + \Delta_r$. 
Moreover, $\msd$ is the corresponding interior lattice point 
in $\intr(\sigma_\Delta)$.  There exists a primitive lattice vector $z'$ in 
the dual lattice $N'$ such that $z' \in \sigma_\Delta^\dual$, $(z')^\perp \cap 
\sigma_\Delta = \Gamma \cap   \sigma_\Delta $, and $\pro{\msd}{z'} =1$. 
We identify the dual lattice $N'$ with $\Nq/\sum_{i,j}
\Z(e_i^* - e_j^*)$. Since $z$ and $z'$ have the same 
orthogonal subspace in $\Mq'$ there exists a positive integer $k$ such that 
$z = kz'$ modulo 
the sublattice generated by $e_2^*-e_1^*, \ldots, e_r^* - e_1^*$.
In this case, we have $\pro{\msd}{z} =k \pro{\msd}{z'} =k$, 
because $\pro{\msd}{ e_i^* - e_j^*} =1 -1 =0$ for all $i,j$. 
Now it remains to show that $k =1$.  
Using the splittings 
\[ \Nq = N \oplus \bigoplus_{i =1}^r \Z e_i^*, \;\;\; 
\bigoplus_{i=1}^r \Z e_i^* = \Z e_1^* \oplus \bigoplus_{i =2}^r 
\Z (e_i^*-e_1^*),   \]  
we can write  
\[ z = kz' + \sum_{i=1}^r l_i e_i^* \]
for some integers $l_1, \ldots, l_r$. Applying $\pro{p_j}{*}$ to 
both sides of the last equation, we get 
\[ 0 = \pro{p_j}{z} = \pro{p_j}{kz'} + \pro{p_j}{\sum_{i=1}^r l_i e_i^*} = 
k\pro{p_j}{z'} + l_j,\]
i.e., $k$ divides $l_j$ for all $ 1 \leq j \leq r$. Since $z =  
kz' + \sum_{i=1}^r l_i e_i^*$  is  a primitive lattice vector, we get $k =1$. 
\end{proof}

\smallskip \subsection{Special simplices}

Our next purpose is to explain the relations between Cayley cones and 
Cayley polytopes to some results in \cite{Ath05,BR05}. 
It would be more convenient to work with polytopes 
in  the dual space $\NR$. For this reason, we denote these polytopes
by $\nabla$. 

Athanasiadis introduced in \cite{Ath05} the concept of a special simplex 
$S$ of a convex polytope $\nabla$. 
Here we slightly modify this definition  by demanding the vertices of the 
simplex $S$ to be lattice points in $\nabla$.

\begin{definition}{\rm Let $\nabla \subset \NR$ be a $d$-dimensional 
lattice polytope. 
A simplex $S$ spanned by $r$ affinely independent 
lattice points in $\nabla$ is called a 
{\em special $(r-1)$-simplex} of $\nabla$, if 
each facet of $\nabla$ contains exactly $r-1$ vertices of $S$. 
}
\end{definition}

Now, we can give an addendum to Cor. \ref{compl}:

\begin{proposition} Let $\nabla \subset \NR$ be a 
Gorenstein polytope of index $r$. There is a $1$-to-$1$ correspondence 
between Cayley polytope structures of length $r$ of the dual Gorenstein 
polytope $\nabla^* \subset \MR$  and special  $(r-1)$-simplices $S 
\subset \nabla$.
A special $(r-1)$-simplex $S \subset \nabla$ 
can be characterized as a lattice $(r-1)$-dimensional simplex in $\nabla$, 
which is not contained in the boundary of $\nabla$. 
Moreover, all special $(r-1)$-simplices in $\nabla$ have the 
same barycenter.
\label{ath}
\end{proposition}

\begin{proof}

Let $\sigma \subseteq \MR$ be the reflexive Gorenstein cone 
associated to $\Delta := \nabla^*$. We identify $\Delta$ with  the support 
 $\Dt$ of $\sigma$ and $\nabla$ with  the support $\nt$ 
of $\sigma^\dual$. If $\Delta$ is a Cayley polytope of length $r$, then, by  
Cor. \ref{compl}, there exist  $n_1, \ldots, n_r \in \nt \cap \Nq$ 
such that $n_1 + \cdots + n_r = \ns$. If  $F$ is  a facet of 
$\nt$, then there is a vertex $v \in \Dt$ such that  $\pro{v}{F} = 0$. 
Let $\Dt = \conv(\Dt_1, \ldots, \Dt_r)$ as in Prop. \ref{char}. 
Then $v \in \V(\Dt_i)$ for some $i$. This implies $\pro{v}{n_i} = 1$ 
and $\pro{v}{n_j} = 0 $ for $j \not= i$. So  
$n_i \not\in F$, but $n_j \in F$ for $j \not= i$. 
So  $\conv(n_1, \ldots, n_r)$ is a special $(r-1)$-simplex in  
$\nabla$.

Now assume that $n_1, \ldots, n_r \in \nt \cap \Nq$ are vertices of a 
 $(r-1)$-dimensional simplex that  is not contained in the boundary 
of $\nabla$. Then  $n_1 + \cdots + n_r \in r \nt$ is 
in the interior of $\sigma^\dual$. Since  $r \nt$ is reflexive, we have 
 $\intr(r \nt) \cap \Nq = \{\ns\}$. Therefore $n_1 + \cdots + n_r = \ns$.
By Prop. \ref{char}, $\Delta$ is a Cayley polytope of length $r$. Moreover,
above arguments show that $S= \conv(n_1, \ldots, n_r)$ is a special 
$(r-1)$-simplex in $\nabla$. If  $S'= \conv(n_1', \ldots, n_r')$ is 
another special simplex in $\nabla$ then 
\[  n_1 + \cdots + n_r = n_1' + \cdots + n_r' = \ns.\]
Thus $S$ and $S'$ have the same barycenter.
\end{proof}

\begin{example} 
{\rm If $\nabla \subset \NR$ is a reflexive polytope 
(Gorenstein polytope of index 
$r =1$) then there exists only one special $(r-1)$-simplex in $\nabla$ which 
is exactly the unique interior lattice point in $\nabla$.} 
\label{r=1}
\end{example}

\begin{example} {\rm Let $\nabla_d$ be the standard 
$d$-dimensional unit cube in $\R^d$:
$\nabla_d := \{ x \in \R^d\; : \; 0 \leq x_i \leq 1 \}$. 
Then $\nabla_d$ is a Gorenstein
polytope of index $r =2$. There exist exactly $2^{d-1}$ different 
special $1$-simplices in $\nabla_d$ which are spanned by all possible 
pairs of opposite 
vertices of $\nabla_d$. The dual Gorenstein polytope $\nabla_d^*$ has exactly 
$d$ special $1$-simplices which $1$-to-$1$ correspond to 
$d$ different Cayley structures 
of the $d$-dimensional unit cube $\nabla_d$. }  
\label{cube}
\end{example}

Let us recall the following definition:

\begin{definition}{\rm 
A lattice polytope $\nabla \subseteq \NR$ is called 
{\em integrally closed}, if 
any lattice point in the Gorenstein cone $\sigma$ over $\nabla$ 
is a  sum of lattice points from the support $\sigma_{(1)} \cong \nabla$.
}
\end{definition}

We get the following corollary:

\begin{corollary}
Let $\nabla \subseteq \NR$ be an integrally closed Gorenstein 
polytope of index $r$. 
Then $\nabla$ contains a special $(r-1)$-simplex and 
$\nabla^*$ is a Cayley polytope of length $r$.
\label{normal}
\end{corollary}

\begin{proof} It follows immediately from the fact that the lattice 
point $\ns$ is a sum of $r$ lattice points $n_1, \ldots, n_r$ from $\nt$. 
\end{proof}

\begin{remark} 
{\rm One can easily show that if a lattice polytope $\nabla$ 
admits a unimodular triangulation, then $\nabla$  is integrally closed.   
Under the assumption that $\nabla$  admits 
a unimodular triangulation the statement in \ref{normal} 
was proved  in \cite[Prop. 3]{Sti98}. 
} 
\end{remark} 
\smallskip

One of the main results in \cite{BR05} is Corollary 7, 
where Bruns and Roemer show that for any integrally closed 
Gorenstein polytope there 
exists a reflexive polytope with the same $h^*$-vector. 
Here we give a slightly more general version of this result 
that does not need any commutative algebra. 
We first recall the definition of the $h^*$-polynomial 
(also called $\delta$-polynomial), cf. \cite{BN06,Hib92,Sta87}:

\begin{definition}{\rm Let $\nabla \subset \NR$ be a $d$-dimensional 
lattice polytope. 
The {\em $h^*$-polynomial} $h^*_\nabla(t)$ of $\nabla$ is defined as the 
numerator of the Ehrhart series, i.e., 
$h^*_\nabla(t) = (1-t)^{d+1} \sum_{k \in \N} \card{k \nabla \cap N} t^k$. 
\label{h-po}
}
\end{definition}

\begin{remark}
{\rm 
The $h^*$-polynomial has degree at most $d$. 
It is well-known, see \cite{Hib92}, 
that $\nabla$ is Gorenstein of index $r$ 
if and only if for the $h^*$-polynomial one has:
\[ t^{d-r+1} h^*_\nabla(t^{-1}) = h^*_\nabla(t).\]
 In this case, $\nabla$ is reflexive if and only if 
the degree of $h^*_\nabla$ is  $d$.} 
\label{symm}
\end{remark}

The following theorem is very close to results in \cite{Ath05}, \cite{BR05}, 
and \cite{OH05}: 

\begin{theorem}
Let $\nabla \subseteq \NR$ be a Gorenstein polytope of index $r$. 
Assume that  $\nabla$ contains a special $(r-1)$-simplex with vertices 
$n_1, \ldots, n_r$, i.e. the dual Gorenstein polytope $\nabla^*$ 
is a Cayley polytope associated with some 
lattice polytopes $\Delta_1, \ldots,\Delta_r$.  
Then the projection of  $\nabla$ along 
$\aff(n_1, \ldots, n_r)$ is  a $(d-r+1)$-dimensional reflexive 
polytope $P$ such that 
$P= (\Delta_1 + \cdots + \Delta_r)^*$ and 
\[h^*_P(t) = h^*_\nabla(t).\]
\label{brtheo}
\end{theorem}

\begin{proof}
The statement follows straightforward from \cite[Lemma 4, Lemma 5]{BR05}. 
Another variant of the same ideas is contained in 
\cite[Prop.3.12]{RW05} and \cite[Lemma 3.4, Cor. 4.1, Cor. 4.2]{Ath05}. 
For convenience, we explain the main combinatorial 
idea of the proof. 

Let $\nt$ be the support of the Gorenstein cone over 
$\nabla$. We identify $\nabla$ with $\nt$.
For any $i \in \{ 1, \ldots ,r\}$ let $\F_i$ be the set 
of facets of $\nt$ that do {\em not} contain $n_i$. 
We denote by $\Gamma$ the union of all faces 
$F_1 \cap \cdots \cap F_r \subset \nt$, where 
$F_i \in \F_i$ for $i = 1, \ldots, r$. 

Take one such a face  
$G := F_1 \cap \ldots \cap F_r \in \Gamma$. 
We define a lattice polytope $G' := \conv(G, n_1, \ldots, n_r)$. 
Let $v_1, \ldots, v_r \in \nt^*$ be the 
vertices corresponding to the facets $F_1, \ldots, F_r$. By Prop. \ref{char}, 
we have $\pro{v_j}{n_i} = \delta_{i,j}$ for $i,j \in \{1, \ldots, r\}$. 
Therefore 
$G'$ is an $r$-fold pyramid over $G$, 
i.e., $h^*_{G'}(t) = h^*_G(t)$, cf. \cite[Remark 2.6]{Bat07}. 
By \cite[Lemma 4]{BR05} the set of all $G'$ (together with  
their subfaces) form a polyhedral decomposition of $\nt$.

Denote by $\pi$  the projection of $\nabla$ along 
$\aff(n_1, \ldots, n_r)$. Let $x \in \intr(\nabla)$ be the 
unique interior point with lattice 
distance $1/r$ from any facet of $\nabla$. 
We define the refined lattice $N' \supseteq N$ with $r = 
[N':N]$ by $N' := N +\Z x$. By  
Prop. \ref{ath} and Prop. \ref{char}, we have  
$(n_1 + \cdots + n_r)/r = x$. So the projection $\pi$ can 
be seen as the surjective map $N' \to N'/(\Z x + \sum_{i=1}^{r-1} \Z n_i)$
where  $x,n_1, \ldots, n_{r-1}$ form a part of $\Z$-basis of $N'$.
Let $F$ be a face of $\nt$ which is the preimage of a facet of $P$. 
By definition of special $(r-1)$-simplex, any facet of $\nt$ is contained in 
$\F_1 \cup \cdots \cup \F_r$. Since $n_i \not\in F$ for 
all $i = 1, \ldots, r$ $(\pi(n_i) =0 \not\in \partial P)$, 
there exists $F_i \in \F_i$ which contains $F$. 
Thus  $F \subseteq G= F_1 \cap \cdots \cap F_r$  for some 
$G \in \Gamma$. On the other hand, $\dim F \geq d-r$ and $\dim G = d-r$. 
Therefore,  $F = G$. By \cite[Lemma 5]{BR05}, the projection 
$\pi\; : \; \nt  \to P$ respects two polyhedra decompositions
\[ \nt = \bigcup_{G \in \Gamma} \conv(G, n_1, \ldots, n_r) \]
\[ P =  \bigcup_{G \in \Gamma} \conv(\pi(G), 0), \] 
where  each $\conv(G, n_1, \ldots, n_r)$ is a  $r$-fold pyramid 
over a face $G \in \Gamma$ and the polytope
 $\conv(\pi(G), 0)$ is $1$-fold pyramid with 
vertex $0$ 
over a facet in $\pi(G) \subset \rand P$. Since  the $h^*$-polynomials of 
 $\conv(G, n_1, \ldots, n_r)$ and $\conv(\pi(G), 0)$ are the same, a standard inclusion-exclusion 
argument shows that the $h^*$-polynomials of $\nt$ and 
$P$ are the same. Using \ref{symm}, this implies that $P$ is a reflexive
polytope. By 
Theorem \ref{sum},  
$\Delta_1 + \cdots + \Delta_r$ is a reflexive polytope. 
Moreover, the proof of Theorem \ref{sum} shows 
that $P$ is the dual reflexive polytope to $\Delta_1 + \cdots + \Delta_r$.
\end{proof}

\begin{example}{\rm 
The above theorem is  illustrated in the following example: 
Consider one-dimensional lattice polytopes   
$\Delta_1:= \conv((-1,0),(0,-1))$ and $\Delta_2:= 
\conv((0,0),(1,1))$. We remark that  
$\Delta_2$ contains $0$, but $\Delta_1$ does not. The Cayley polytope 
$\Dt = \Delta_1*\Delta_2$ is a Gorenstein polytope of index $2$. 
On the other hand, $\Dt^*$ is isomorphic to $2 S_3$, where $S_3$ 
is the unimodular $3$-simplex. The polytope  $\Dt^*$ contains 
a special $1$-simplex $S$ corresponding to the Cayley polytope structure 
of $\Dt$ (the two vertices of $S$ are 
marked by  boxes). By projecting along this simplex, we obtain 
the reflexive square $P$. The boundary subcomplex $\Gamma$ of 
$\Dt^*$ (which maps  bijectively  to the 
boundary of $P$) is marked by  thick edges. 
The dual reflexive polytope $P^*$ of $P$ is $\Delta_1 + \Delta_2$. It is 
the intersection of $2\Dt$ with an affine plane which is orthogonal 
to the special simplex $S$.  \\

{\centerline{ 
\psfrag{a}{$\Delta_2$}
\psfrag{b}{$\Delta_1$}
\psfrag{c}{$\Dt = \Delta_1 * \Delta_2$}
\psfrag{d}{$\Delta_1 + \Delta_2$}
\psfrag{e}{$2 \Dt$}
\psfrag{f}{$\Dt^*$}
\psfrag{g}{$P$}
\psfrag{h}{$2 \Dt^*$}
\psfrag{i}{$\twoheadleftarrow$}
\psfrag{j}{$\hookrightarrow$}
\psfrag{k}{$\twoheadrightarrow$}
\includegraphics{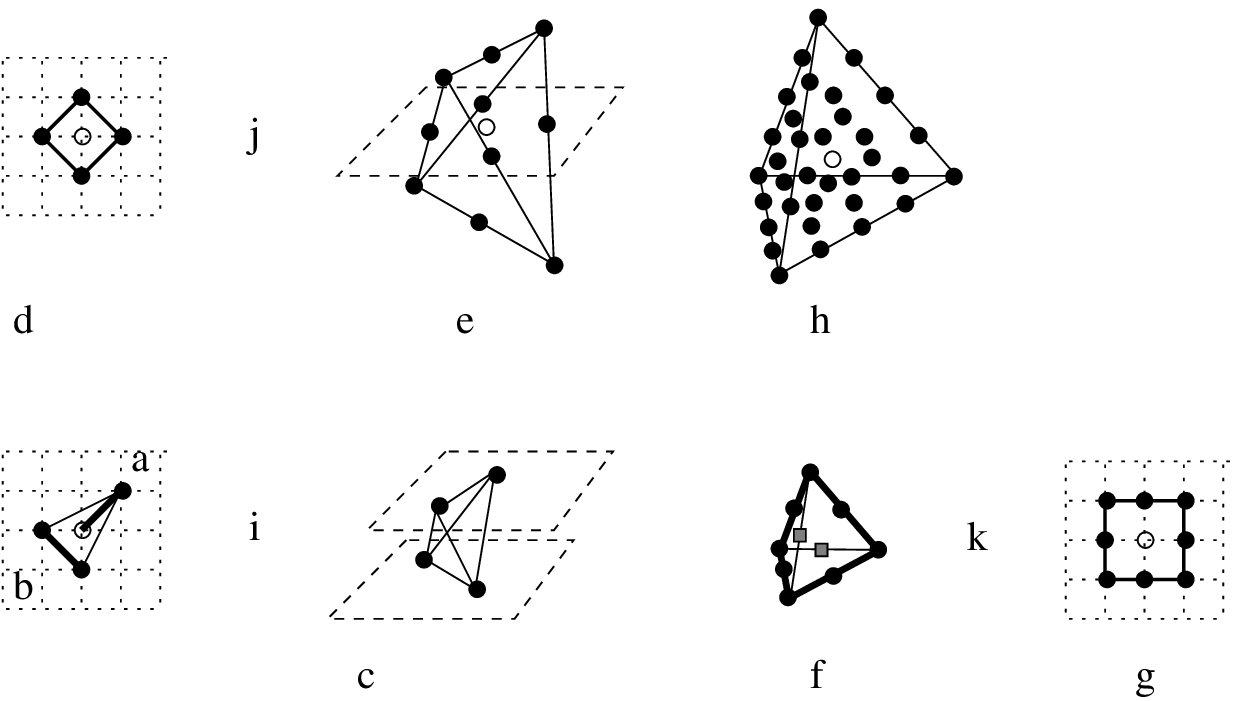}}}
\label{fig1}
}
\end{example}

\section{Combinatorics of  nef-partitions}

\subsection{Characterization of nef-partitions}

Now we are interested in the case when  both dual reflexive Gorenstein cones 
are completely split. Using special simplices, we want to obtain 
a simple combinatorial characterization of nef-partitions and their 
duality.

\smallskip

Let us consider the following purely combinatorial definition of 
a nef-partition (cf. \cite[Prop.3.2]{KRS03}):

\begin{definition}{\rm 
Let $\Delta \subset \MR $ be a $d$-dimensional reflexive polytope and  
$m = \intr(\Delta) \cap M$. A Minkowski sum decomposition 
$\Delta = \Delta_1 + \cdots + \Delta_r$ where $\Delta_1, \ldots, \Delta_r 
\subseteq \MR$ are lattice polytopes   is called a {\em nef-partition} 
of $\Delta$ of 
length $r$, if there exist lattice points $p_1 \in \Delta_1, \ldots, p_r 
\in \Delta_r$ such that $p_1 + \cdots + p_r = m$. Since the reflexive 
polytope $\Delta$ is determined by $\Delta_1, \ldots, \Delta_r$, for 
simplicity, we the name {\em nef-partition} just for the the set 
of lattice polytopes $\Delta_1, \ldots, \Delta_r$.   

We call a nef-partition {\em centered}, if $p_i =0$ for all $1 \leq i \leq r$.
We call a nef-partition {\em proper}, if $\dim \Delta_i > 0$ 
for all $1 \leq i \leq r$. 
\label{nefdef}
}
\end{definition}

\begin{example} 
{\rm Let $\Delta \subset \R^d$ be the $d$-dimensional reflexive cube defined by 
the conditions $|x_i| \leq 1$ $( 1 \leq i \leq d)$. Then $\Delta$ is a 
Minkowski sum of two unit cubes $\Delta_1$ and $\Delta_2$ defined as 
\[ \Delta_1 : = \{ x=(x_i) \in \R^d\; : \; 0 \leq x_i \leq 1 \; 
\forall i \}, \] 
\[ \Delta_2 : = \{ x=(x_i) \in \R^d\; : \; -1 \leq x_i \leq 0 \; 
\forall i \}. \] 
Since $0 \in \Delta_i$ $( i =1,2)$, the Minkowski sum  $\Delta = \Delta_1 + 
\Delta_2$ defines   a proper centered nef-partition of length $2$. This 
example comes from the Gorenstein polytopes considered in \ref{cube}.} 
\label{nefp}
\end{example}  

\begin{remark} 
{\rm It clear that if $\Delta = \Delta_1 + \cdots + \Delta_r$ 
is a nef-partition of $\Delta$, then 
\[ (\Delta -m) = (\Delta_1 -p_1) + \cdots + (\Delta_r -p_r) \]
is a centered nef-partition. Moreover, if 
\[ \Delta = \sum_{i=1}^r \Delta_i \]
is a centered nef-partition of $\Delta$, then 
\[ \Delta = \sum_{i\; : \; \Delta_i \neq 0} \Delta_i \]
is a proper centered nef-partition.
Therefore, the study of arbitrary nef-partitions immediately reduces to 
study of proper centered nef-partitions. } 
\end{remark} 

\begin{remark}
{\rm 
The above definition of proper centered 
nef-partitions is nearly precisely the one given by Kreuzer, 
Riegler and Sahakyan
in \cite[Prop.3.2]{KRS03}. They also demand an  additional assumption 
$\Delta_i \cap \Delta_j = \{0\}$ for $i \not= j$, 
which is actually unnecessary (see Cor.~\ref{inter}).
} 
\end{remark}

\begin{remark} 
{\rm  From the viewpoint of algebraic geometry, proper 
centered nef-partitions were introduced by Borisov in \cite{Bor93}. 
Let  $\P_{\Delta}$ be  a Gorenstein toric Fano variety corresponding 
to the reflexive polytope $\Delta$ which contains $0$ in its 
interior. Denote by $\T$ the $d$-dimensional 
torus acting on  $\P_{\Delta}$ (we identify  $\T$ with 
 its open dense orbit in  $\P_{\Delta}$).  A Minkowski sum 
decomposition $\Delta = \Delta_1 + \cdots + \Delta_r$ 
describes the anticanonical sheaf  ${\mathcal O}(-K_{\P_\Delta})$ as a 
tensor product  of $r$ $\T$-equivariant semi-ample  
line bundles  ${\mathcal L_1}, \ldots, 
{\mathcal L_r}$. In this case, the lattice points $p_i \in \Delta_i$ can be 
identified with some global $\T$-invariant 
sections $s_i$ of ${\mathcal L}_i$ which do not 
vanish on $\T \subset \P_{\Delta}$. The condition 
$\protect{p_1 + \cdots + p_r =0}$ 
means that if $D_i$ is the zero set of the global section $s_i$ $(i =1,
\ldots, r)$, then 
$D_1 + \cdots + D_r$ is the anticanonical divisor of $\P_{\Delta}$ 
which contains every irreducible component of $\P_{\Delta} \setminus 
\T$ with multiplicity $1$. 
In particular, $D_i$ and $D_j$ do not have 
common irreducible components for $i \neq j$. Moreover, $\dim\, \Delta_i 
\neq 0$ if and only if $D_i \neq 0$. 
Therefore, if $\dim\, \Delta_i 
\neq 0$ for all $1 \leq i \leq r$, then  we get a {\em partition} 
of the set of irreducible components of $\P_{\Delta} \setminus 
\T$ into $r$ nonempty subsets which 
define {\em nef-divisors} $D_1, \ldots, D_r$ such 
that $D_1 + \cdots + D_r$ is the anticanonical divisor. This explains 
the origin of the notion {\em nef-partition}   in \cite{Bor93}. 
 } 
\end{remark}

We characterize nef-partitions in the following way:

\begin{proposition}
Let $\Delta$ be a reflexive polytope with the unique interior 
lattice point $m$, and let $\Delta_1, \ldots, \Delta_r$ be lattice 
polytopes such that $\Delta = \Delta_1 + \cdots + \Delta_r$. 
Denote by  $\sigma$ the Cayley cone associated to 
$\Delta_1, \ldots, \Delta_r$. Then the following statements are equivalent:
\begin{enumerate}
\item the dual cone $\sigma^\dual$ is a completely 
split reflexive Gorenstein 
cone of index $r$
\item $\Delta_1 * \cdots * \Delta_r$ is a Gorenstein polytope 
of index $r$ containing a special $(r-1)$-simplex
\item there exist lattice points 
$p_i \in \Delta_i$ $(1 \leq i \leq r)$ such that $p_1 + \cdots + p_r = m$, 
i.e., $\Delta =  \Delta_1 + \cdots + \Delta_r$ is a nef-partition.
\end{enumerate}
\label{splittheo}
\end{proposition}

\begin{proof}

(1) $\ra$ (2): If $\sigma^\dual$ is a reflexive Gorenstein cone 
of index $r$, then this holds also for $\sigma$. Hence Theorem \ref{sum} 
implies that $\Delta_1 * \cdots * \Delta_r$ is a Gorenstein
 polytope of index $r$. Since $\sigma^\dual$ is completely split, 
Prop. \ref{ath} implies that $\Delta_1 * \cdots * \Delta_r$ 
contains a special $(r-1)$-simplex.

(2) $\ra$ (3):  Since 
$\sigma^\dual$ is completely split, 
Cor. \ref{compl} implies that there are $\p_1, \ldots, \p_r \in \Dt \cap 
\Mq$ 
such that $\p_1 + \cdots + \p_r = \msd$. Let $e^*_1, \ldots, e^*_r \in \Dt^* \cap \Nq$ as in Prop. \ref{char} and 
Cor. \ref{compl}. Since $\pro{\msd}{e_i^*} =1$ 
$\forall i$ and $\pro{\p_j}{e_i^*}$ is a nonnegative integer $\forall j$, 
every polytope $\Dt_i = \Delta_i \times e_i$ contains exactly one point from 
$\{\p_1, \ldots, \p_r\}$. 
Without loss of generality we can assume that $\p_i \in \Dt_i$ $\forall i$.
We can write $\p_i = p_{i} \times 
e_{i} \in \Dt_{i} \cap \Mt$ $( 1\leq i \leq r)$.  
By \ref{sum},   $\p_1 + 
\cdots + \p_r = \msd = m \times (e_1 + \cdots + e_r)$. 
Hence for $i=1, \ldots, r$ we have $p_{i} \in \Delta_{i} 
\cap M$ with $p_1 + \cdots + p_r = m$ and  $\Delta_1, \ldots, \Delta_r$ is 
a nef-partition.

(3) $\ra$ (1): By Theorem \ref{sum} we know that $\sigma$ is a 
reflexive Gorenstein cone of index $r$, hence also $\sigma^\dual$ is. 
We define $\p_i := p_i \times e_i \in \Mq = M \oplus 
\Z^r$ for $i=1, \ldots, r$. Then 
$\p_1, \ldots, \p_r \in \Delta_1 * \cdots * \Delta_r \cap \Mq$ with 
$\p_1 + \cdots + \p_r = \msd$. Hence $\sigma^\dual$ is completely
 split by Cor. \ref{compl}.
\end{proof}

We can sum this up in a more symmetric way:

\begin{corollary}
A reflexive Gorenstein cone $\sigma$ is associated to a 
nef-partition if and only if both cones 
 $\sigma$ and $\sigma^\dual$ are completely split.

A Gorenstein polytope $\Dt$ of index $r$ is the Cayley polytope of a 
nef-partition if and only if $\Dt$ and $\Dt^*$ are Cayley polytopes 
of length $r$, 
or equivalently, if and only if both polytopes 
$\Dt$ and $\Dt^*$ contain special $(r-1)$-simplices.
\label{splitcoro}
\end{corollary}

\begin{example}{\rm
In Example \ref{fig1}, the polytope  $\Dt = \Delta_1 * \Delta_2$ 
is a Gorenstein polytope of index $2$. It does not contain a 
special $1$-simplex. 
Thus $\Dt^*$ is not a Cayley polytope of length $2$. 
In other words, the Cayley cone $\sigma$ associated to $\Delta_1, \Delta_2$ 
is completely split, while $\sigma^\dual$ is not. Hence the sum 
$\Delta_1 + \Delta_2$ is not a nef-partition of the reflexive polytope 
$\Delta_1 + \Delta_2$, because  
there are no lattice points $p_1 \in \Delta_1$, $p_2 \in \Delta_2$ such that 
 $p_1 + p_2 = 0$.
}
\end{example}

From Cor. \ref{normal} we get:

\begin{corollary}
Let $\Dt$ be a Gorenstein polytope such that both 
$\Dt$ and $\Dt^*$ are  integrally closed. Then $\Dt$ (and also $\Dt^*$) 
is a Cayley polytope associated to 
a nef-partition.
\end{corollary}

\begin{remark}{\rm 
In dimension $2$, all lattice polygons 
$\Delta_1, \Delta_2$ such that $\Delta_1 + \Delta_2$ is reflexive 
were listed in \cite[Fig.4,5]{ET05}. 
In dimension $3$ and $4$, all reflexive polytopes are 
completely classified in \cite{KS98,KS00}. Using the software 
package PALP, one can find all possible 
centered nef-partitions whose Minkowski sum yields a given 
reflexive polytope \cite{KS04}. However, it would be desirable to have an 
efficient algorithm that classifies all centered nef-partitions 
of given length $r$ in fixed dimension $d$. At the moment, the only way is 
to compute {\em all} reflexive polytopes of dimension $d$, 
which is practically impossible for $d \geq 5$.
}
\end{remark}

\smallskip \subsection{Duality of  nef-partitions}

\begin{definition}
{\rm 
Let $P$ and $P^*$ be two dual to each other $d$-dimensional 
Gorenstein polytopes
of index $r$. Assume that there exist two special $\protect{(r-1)}$-simplices 
$S \subset P$ and $S' \subset P^*$. By \ref{ath}, these 
simplices define Cayley structures on polytopes $P$ and $P^*$. 
Therefore,  for some lattice polytopes $\Delta_1, \ldots, \Delta_r$ and 
$\nabla_1, \ldots, \nabla_r$  we obtain 
\[ P = \Delta_1 * \cdots * \Delta_r, \;\; P^* = \nabla_1 * \cdots *\nabla_r, \]
and two nef-partitions 
\[ \Delta = \Delta_1 + \cdots + \Delta_r, \;\; \nabla = \nabla_1 + 
\cdots + \nabla_r \]
of $(d-r+1)$-dimensional reflexive polytopes $\Delta$ and $\nabla$. 
We call these two nef-partitions {\em dual to each other}.  
}
\label{nef-part-dual} 
\end{definition}

\begin{remark} 
{\rm The duality of nef-partitions defined in \ref{nef-part-dual} looks
trivial, because it simply exchanges Gorenstein polytopes 
 $P$ with  $P^*$ and special $(r-1)$-simplices $S$ and $S'$. By \ref{r=1},
for $r=1$ this duality reduces to the polar duality between reflexive 
polytopes $P$ and $P^*$. However, Definition  \ref{nef-part-dual} is not
very convenient if we want to determine explicitly 
the polytopes $\nabla_1, 
\ldots, \nabla_r$ from polytopes $\Delta_1, \ldots, \Delta_r$. 
For this reason, we consider below another combinatorial approach to  
the duality of centered  nef-partitions.}
\end{remark}   

\begin{proposition}
Let $\Delta \subset \MR$ be a reflexive 
polytope and $0 = \intr(\Delta) \cap M$. 
Assume that  $\Delta= \Delta_1 + \ldots +  \Delta_r$ is a centered 
nef-partition. 
Then  $r$ polytopes
\[\nabla_i := \{y \in \NR \::\: \pro{\Delta_j}{y} \geq -\delta_{i j} 
\;\forall\, j=1, \ldots, r\}, \;\; i\in \{ 1, \ldots, r\}.\]
define the dual centered nef-partition.
\label{nefdualdef}
\end{proposition}

\begin{proof}  
Consider two 
lattices 
$\Mq = M \oplus \Z e_1 \oplus \cdots \oplus \Z e_r$, and 
$\Nq = N \oplus \Z e^*_1 \oplus \cdots \oplus \Z e^*_r$, where 
$\pro{e_i}{e^*_j} = \delta_{i j}$.
Let 
$\sigma \subseteq \MRq$ be the reflexive Cayley cone associated with the 
lattice polytopes $\Delta_1, \ldots, \Delta_r \subseteq \MR$. Denote by 
$\sigma^\dual \subset \NRq$ the dual reflexive cone. We set 
$\Dt_i = \Delta_i \times e_i$. Then 
 $0 \times (e_1 + \cdots + e_r) = \msd$, and 
$0 \times (e^*_1 + \cdots + e^*_r) = \ns$. Let $\Dt^*$ be 
the support of $\sigma^\dual$. 

For every $i = 1, \ldots, r$,  we define polytopes 
$\nt_i := \{y \in \Dt^* \;:\; 
\pro{e_j}{y} = 0$ for $j \not= i$$\}$.  
By Prop. \ref{char}, we have  
$\Dt^* = \conv(\nt_1, \ldots, \nt_r)$. 
By Prop. \ref{splittheo}, for   every $i = 1, \ldots, r$, we can write 
$\nt_i$ as $\nabla_i '\times e^*_i$ for some  $\nabla_i' \subseteq \NR$. 
We want  to show that $\nabla_i = \nabla_i'$ $( 1\leq i \leq r)$. This 
follows  from the observation that  
for  $y \in \NR$ and $i \in \{1, \ldots, r\}$ one has  
\[\pro{\Delta_j}{y} \geq -\delta_{i j} \;\forall\, j=1, \ldots, r\]
if and only if
\[\pro{\Dt_j}{y \times e^*_i} \geq 0 \;\forall\, j=1, \ldots, r\]
(i.e.,  $y \times e^*_i \in \nt_i$ $\Leftrightarrow$ $y \in \nabla_i$).
Obviously, one has $0 \in \nabla_i$ $(1 \leq i \leq r)$. So we obtain that 
$\nabla_1, \ldots, \nabla_r \subseteq \NR$ is the centered 
dual nef-partition. 
\end{proof}

\begin{example}
\label{dual-nefp} 
{\rm The dual centered nef-partition to the one 
considered in \ref{nefp} consists 
of two $d$-dimensional standard simplices 
\[ \nabla_1 := \{ x = (x_i) \in \R^d\; : \; x_1 + \ldots + x_d \geq -1,\, 
x_i \leq 0\; \forall i \}, \]
\[ \nabla_2 := \{ x = (x_i) \in \R^d\; : \; x_1 + \ldots + x_d \leq 1,\, 
x_i \geq 0\; \forall i \}. \]} 
\end{example} 

\begin{example}
{\rm The following picture below  
describes another example of two dual to each other nef-partitions 
in dimension $2$. \\
{\centerline{
\psfrag{a}{$\Delta_2$}
\psfrag{b}{$\Delta_1$}
\psfrag{c}{$\nabla_1$}
\psfrag{d}{$\nabla_2$}
\includegraphics{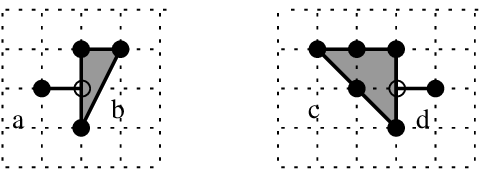}}}
 In this case, the dimensions of $\Delta_2$ and $\nabla_2$ 
are not maximal. 
  \label{fig3}
\label{nef-p2}
}
\end{example}

\smallskip

For the remainder of this section let $\Delta_1, \ldots, 
\Delta_r$ and $\nabla_1, \ldots, \nabla_r$ 
be two dual to each other centered nef-partitions.

\smallskip

\begin{corollary}
For any $i \in \{ 1, \ldots, r\}$ one has $\Delta_i =0$ if and only if 
$\nabla_i = 0$. In particular, 
the dual nef-partition of a proper nef-partition is a proper nef-partition.
\label{proper}
\end{corollary}

\begin{proof}
Without loss of 
generality, we can assume that $\Delta_1 =0$. Then the Minkowski sum 
$\Delta_2 + \cdots + \Delta_r = \Delta_1 + \cdots + \Delta_r$ is a 
reflexive polytope $\Delta$ containing $0$ in its interior.     
Hence $\{ y \in \NR\; : \; \pro{x}{y} \geq 0 \; \forall x \in \Delta \} = 0$, 
so 
by the formulas in \ref{nefdualdef} we obtain that $\nabla_1 = 0$.  
By duality of nef-partitions, the condition  
$\nabla_1 = 0$ implies $\Delta_1 =0$, analogously. 
\end{proof}

With the same argument we also derive from Prop. \ref{nefdualdef}:

\begin{corollary}
$\pos(\nabla_i) \cap \pos(\nabla_j) = \{0\}$ for $i \not= j$.
\label{inter}
\end{corollary}

The crucial property of the duality of nef-partitions is that it exchanges 
Minkowski sum and convex hull:

\begin{proposition}
$\conv(\Delta_1, \ldots, \Delta_r)$ is a reflexive polytope which is 
dual to $\nabla_1 + \cdots + \nabla_r$. Moreover, the two lattice polytopes 
$\Delta_1 * \cdots * \Delta_r$ and 
$\conv(\Delta_1, \ldots, \Delta_r)$ have the same $h^*$-polynomial.
\label{minkdual}
\end{proposition}

\begin{proof} 
The statements follow immediately from Theorem \ref{brtheo}, because
\linebreak $\conv(\Delta_1, \ldots, \Delta_r)$ is the projection of 
$\Delta_1 * \cdots * \Delta_r$ along $\aff(e_1, \ldots, e_r)$. 
\end{proof}

Hence, we have for $\Delta = \Delta_1 + \cdots + \Delta_r$ and 
$\nabla = \nabla_1 + \cdots + \nabla_r$:
\[\Delta^* = \conv(\nabla_1, \ldots, \nabla_r) \text{ and } 
\nabla^* = \conv(\Delta_1, \ldots, \Delta_r).\]
\pagebreak

Now, we show that in order to determine the dual nef-partition it suffices 
to dualize the Minkowski sum $\Delta$, and to choose 
the vertices of $\nabla_1, \ldots, \nabla_r$ among the vertices of $\Delta^*$. 
For this we take a closer look at boundary of the lattice polytopes 
$\nabla_1, \ldots, \nabla_r$:

\begin{proposition}
We have for $i = 1, \ldots, r$:
\[\nabla_i \cap \rand \Delta^* = \{v \in \Delta^* \,:\, 
\min_{u \in \Delta_i} \pro{u}{v} = -1\}.\]
Any face $F$ of $\nabla_i$, with $0 \not\in F$, is a face of 
$\Delta^* = \conv(\nabla_1, \ldots, \nabla_r)$ and contained in 
a face of $\nabla = \nabla_1 + \cdots + \nabla_r$. 
\label{faceprop}
\end{proposition}

\begin{proof}

Let $v \in \rand \Delta^*$. Then 
\[(\min_{u \in \Delta_1} \pro{u}{v}) + \cdots + 
(\min_{u \in \Delta_r} \pro{u}{v}) = \min_{u \in \Delta} \,\pro{u}{v} = -1, \]
and these are all numbers in $[0,-1]$. From Prop. 
\ref{nefdualdef} the equation follows.

Now, let $F$ be a facet of $\nabla_i$, with $0 \not\in F$. 
Let $v$ be a point in the relative interior of $F$. 
By Prop. \ref{nefdualdef} we have $\min_{u \in \Delta_i} \pro{u}{v} = -1$, 
so there is a vertex 
$u$ of $\Delta_i$ such that $\pro{u}{v} = -1$. Therefore $F = \{y_i \in 
\nabla_i \,:\, \pro{u}{y_i} = -1\}$, and $F$ is also contained in 
the face $\{y \in \nabla_1 + \cdots + \nabla_r \,:\, \pro{u}{y} = -1\}$ 
of $\nabla_1 + \cdots + \nabla_r$. 

We show that 
$F = \{y \in \Delta^* \,:\, \pro{u}{y} = -1\}$. Let $y \in \Delta^*$ 
with $\pro{u}{y} = -1$. 
There is a convex combination $y = \sum_{j=1}^r \lambda_j v_j$, 
with $\lambda_j \geq 0$, $\sum_{j=1}^r \lambda_j = 1$, 
for $v_j \in \nabla_j$ ($j = 1, \ldots, r$). Hence from $-1 = 
\pro{u}{y}$ we get $-1 \geq  \lambda_i \pro{u}{v_i} \geq - \lambda_i$, 
because of Prop. \ref{nefdualdef}, thus $\lambda_i = 1$, so $y = v_i 
\in \nabla_i$. This shows $y \in F$.
\end{proof}

\begin{remark}{\rm The reader further interested in faces of 
nef-partitions and their dualities is encouraged to consult 
\cite[Prop.2.3]{HZ05} and 
\cite[2.7-2.9]{Gro05}.
}
\end{remark}

As mentioned above, here is a direct way to determine the 
vertices of the dual nef-partition:

\begin{corollary}
We have for $i = 1, \ldots, r$:
\[\V(\nabla_i)\backslash\{0\} = \{v \in \V(\Delta^*) \;:\; 
\min_{u \in \Delta_i} \pro{u}{v} = -1\}.\]
\label{nabcor}
\end{corollary}

\begin{remark}{\rm The previous statement implies that any 
non-zero vertex of $\nabla_i$ ($i = 1, \ldots, r$) is 
a vertex of $\conv(\nabla_1, \ldots, \nabla_r)$. It is not 
difficult to prove, that this also holds under the weaker assumption that 
$\nabla_1, \ldots, \nabla_r$ are lattice polytopes, each 
containing $0$, such that $\nabla_1 + \cdots + \nabla_r$ 
has no interior lattice point except $0$.
}
\end{remark}

Not only vertices but also lattice points of $\Delta^* = 
\conv(\nabla_1, \ldots, \nabla_r)$ necessarily belong to 
one of the lattice polytopes 
$\nabla_1, \ldots, \nabla_r$:

\begin{corollary}
Any non-zero lattice point of $\Delta^*$ is contained for some 
$i \in \{1, \ldots, r\}$ in a face of $\nabla_i$ that does 
not contain the origin. 
In particular $\card{\Delta^* \cap N} = \card{\nabla_1 \cap N} + 
\cdots + \card{\nabla_r \cap N} - r  + 1$.
\end{corollary}

\begin{proof}

Let $v \in \Delta^* \cap N$, $v \not= 0$. There has to be some 
$u \in \V(\Delta_i)$ ($i \in \{1, \ldots, r\}$) with 
$\pro{u}{v}<0$, thus $\pro{u}{v} = -1$. Now use Prop. \ref{faceprop}.
\end{proof}

\section{The $\Est$-function of a Gorenstein polytope}

Borisov and the first author 
gave in \cite{BB96b} an explicit formula for the generating function 
of stringy Hodge numbers of Calabi-Yau complete intersections 
given by a 
nef-partition (the 
stringy Hodge numbers were introduced in \cite{BD96, Bat98}). 
We remark that this formula can be used for arbitrary reflexive 
Gorenstein cones or for arbitrary Gorenstein polytopes. 
However, it is not clear 
a priori whether the corresponding $\Est$-function is a polynomial. 
We expect that it is always the case.

In this section,  we review the above formula using 
$\St$-polynomials and ideas from \cite{BM03}. 
Our  goals are to introduce the $\Est$-function of a Gorenstein polytope, and 
to show some of its properties, and to formulate related open questions 
and conjectures. 

\smallskip

\subsection{The $\St$-polynomial of a lattice polytope}

We start by recalling the definition of the $g$- and $h$-polynomials, 
defined by Stanley 
for general Eulerian posets in \cite{Sta87} (see also \cite[Def.5.3]{BM03}):

\begin{definition}{\rm 
Let $\calP$ be an Eulerian poset of rank $d$, $\hat{0}$ its  minimal element, 
 $\hat{1}$ its maximal element. 
Then the polynomials $g_\calP(t), h_\calP(t) \in \Z[t]$ are  defined 
recursively as follows.  If  $d=0$ then we set  
$g_\calP(t) := h_\calP(t) := 1$. If  $d > 0$ then 
the $h$-polynomial is defined as 
$$h_\calP(t) = \sum_{ \hat{0} <  x \leq  
\hat{1}}\, (t-1)^{\rk(x)-1} \,g_{[x,\hat{1}]}(t).$$ 
For $d > 0$, we define the $g$-polynomial 
$g_\calP(t) = \sum_i g_i t^i$ as a polynomial of degree $< d/2$ whose 
coefficients $g_i$ 
(for $0 \leq i < d/2$) 
are equal to  the coefficients $c_i$ of 
the polynomial $$(1-t) h_\calP(t) = \sum_i c_i t^i,$$ 
i.e. , $g_i = h_i - h_{i-1}$ 
for $0 \leq i < d/2$.
}
\end{definition}

A recent account on this important combinatorial notion 
can be found in \cite{Bra05}. Here we need only the following 
well-known properties:

\begin{remark}  
{\rm The polynomials $g_\calP(t)$ and $h_\calP(t)$ have the followng 
properties: 
\begin{enumerate}
\item every  coefficient of $g_\calP(t)$ is a non-negative integer;
\item $g_\calP(t) = 1$ if and only if the poset $\calP$ is boolean, i.e., 
it is isomorphic to the poset of all faces (including the empty one) of 
a $(d-1)$-dimensional simplex;
\item the polynomial $g$ is multiplicative, i.e., 
$g_{\calP \times \calP'}(t) = g_\calP(t) \cdot  g_{\calP'}(t)$.
\end{enumerate}
}
\end{remark}

Combinatorial data contained in the $g$-polynomial will be combined with 
the enumeration of lattice points via the $h^*$-polynomial in the following  
definition \cite[Def.5.3]{BM03}:

\begin{definition}{\rm Let $P$ be a lattice polytope. 
We define the $\St$-polynomial:
\[\St(P,t) := \sum_{\emptyset \leq F \leq P} 
(-1)^{\dim(P)-\dim(F)} \; h^*_F(t) \; g_{[F,P]}(t),\]
where 
$h^*_F(t)$ is the $h^*$-polynomial of the lattice polytope 
$F$ (with convention $h^*_\emptyset(t) := 1$, see Def. \ref{h-po}) 
and $g_{[F,P]}(t)$ is the $g$-polynomial of the poset 
$[F,P]$ of faces between $F$ and $P$. Note that 
$\St(\emptyset, t) = 1$. However, $\St(P,t) = 0$, if $\dim(P) = 0$.
}
\end{definition}

\begin{remark}{\rm The $\St$-polynomial of a lattice polytope $P$ has 
the following properties:
\begin{enumerate}
\item all coefficients of $\St(P,t)$ are  non-negative integers 
(this follows from \cite[Prop.5.5]{BM03}, 
where it is shown that $\St(P,t)$ is the Hilbert function 
of a graded vector space). 
\item one has  the  reciprocity-law (see \cite[Remark 5.4]{BM03}):
\[\St(P,t) = t^{\dim(P)+1} \, \St(P,1/t).\]
\end{enumerate}
}
\end{remark}

Here is another nice property of the $\St$-polynomial:

\begin{lemma}
Let $P$ be a lattice polytope. If $P$ is a lattice pyramid, then 
$\St(P,t) = 0$.
\label{pyra}
\end{lemma}

\begin{proof}

Let $P$ be a lattice pyramid over the facet $P'$ with apex
 $v$, i.e., $h^*_{P'}(t) = h^*_P(t)$. 
Any face $F$ of $P$, that is not contained in $P'$, 
is a pyramid over a face $F'$ of $P'$ with apex $v$. Hence, 
$h^*_F(t) = h^*_{F'}(t)$ and 
$g_{[F,P]}(t) = g_{[F',P']}(t)$. On the other hand, if 
$F$ is a face of $P'$, then $g_{[F',P]}(t) = g_{[F',P']}(t)$, since 
$[F',P]$ is a product of 
$[F',P']$ and an interval, and $g$ is multiplicative. This yields:
\[\St(P,t) = \St(P',t) + (-1)^{\dim(P)-\dim(P')} \St(P',t) = 0.\]
\end{proof}

In the case of a simplex, the $\St$-polynomial is easy to calculate:

\begin{proposition}
Let $P \subseteq \MR$ be a lattice simplex. Let $\Pt := P \times 1 
\subseteq \MRq = \MR \oplus \R$. 
We denote by $\Pi$ the parallelepiped spanned by the vertices of $\Pt$. Then
\[\St(P,t) = \sum_{x \in \intr(\Pi) \cap M} t^{x_{d+1}}.\]
The degree of $\St(P,t)$ equals at most the degree of $h^*_P(t)$.
\label{simplexprop}
\end{proposition}

\begin{proof}

For a simplex $\Pt$, the $h^*$-polynomial can be 
calculated by counting the lattice points in the 
parallelepiped $\Pi(\Pt) := \{\sum_{i=0}^d \lambda_i v_i \,:\, 0 
\leq \lambda_i < 1\}$, 
where $v_0, \ldots, v_d$ are the vertices of $\Pt$:
\[h^*_\Pt(t) = \sum_{x \in \Pi(\Pt) \cap M} t^{x_{d+1}}.\]
Since the $g$-polynomial of a simplex equals $1$, we get
\[\St(P,t) = \sum_{\emptyset \leq F \leq \Pt} (-1)^{\dim(\Pt)-\dim(F)} 
\; h^*_F(t).\]
Now, for $x \in \Pi(\Pt)$, we denote by $F_x$ the smallest face of 
$\Pt$ such that $x \in \pos(F_x)$. 
Hence,
\[\St(P,t) = \sum_{x \in \Pi(\Pt) \cap M} 
\left(\sum_{F_x \leq F \leq \Pt} (-1)^{\dim(\Pt)-\dim(F)}\right) t^{x_{d+1}}.\]
The Euler-Poincar\'e-formula implies that the right side vanishes, 
whenever $F_x \not= \Pt$. 
From this the equation follows.
\end{proof}

In general, the degree of the $\St$-polynomial of a lattice simplex 
is not equal to the 
degree of the $h^*$-polynomial. This follows immediately from 
Lemma \ref{pyra}.

\begin{example}
{\rm For any lattice polytope $P \subseteq M_{\R}$ we denote
 by $l^*(P)$ (resp. by $l(P)$) 
the number of $M$-lattice points in the relative interior of 
$P$ (resp. in $P \cap M$). 

Now, let $P$ be a $2$-dimensional lattice polytope with $k$ vertices 
(a lattice $k$-gon). In this case, we calculate from the definition 
of the $h$- and $g$-polynomial, 
$h_{[\emptyset,P]}(t) = k + k (t-1) + (t-1)^2 = 1 + (k-2)t + t^2$, thus 
$g_{[\emptyset,P]}(t) = 
1 + (k-3) t$. Then 
\[ \tilde{S}(P,t) = h^*_P(t) - \left( \sum_{{\rm dim}\, \Gamma =1} h^*
(\Gamma) \right)  + 
k - ( 1 + (k-3)t) = \]
\[ = 1 + (l(P)-3)t + l^*(P)t^2 -  \left(  
\sum_{{\rm dim}\, \Gamma =1} 1 + (l^*(\Gamma)-1)t  \right) + k - 
( 1 + (k-3)t) =\]
\[ = l^*(P)( t + t^2), \]
since 
$$l(P) = l^*(P) + \sum_{{\rm dim}\, \Gamma =1} l^*(\Gamma) + k. $$
More generally, if $P$ is a $d$-dimensional lattice 
polytope and $l^*(P) \neq 0$, then $ \tilde{S}(P,t)$ and $h^*_P(t)$ have 
the same leading term $l^*(P) t^d$, because 
\[ {\rm deg} \,(h^*_{F}(t) \, g_{[F, P]}(t)) \leq d-1 \]
for any proper face $F \subset P$, since 
$\deg(h^*_F(t)) \leq \dim(F)$ and $\deg(g_{[F, P]}(t)) < (\rank [F,P])/2 = 
(d-\dim(F))/2$.
} 
\end{example}    

\smallskip

\subsection{ $\Est$-function  of  a Gorenstein polytope}

\begin{definition}{\rm Let $\Delta$ be a $d$-dimensional 
Gorenstein polytope of index $r$. 
We define the following function of two independent variables 
$u$ and $v$: 
\[\Est(\Delta;u,v) := (uv)^{-r} \sum_{\emptyset \leq F \leq \Delta} 
(-u)^{\dim(F)+1} \; \St(F,u^{-1} v) \; \St(F^*,u v),\]
where $F^*$ is the face of $\Delta^*$ dual to $F$.
We call it the {\em stringy $E$-function} (or  $\Est$-function) of $\Delta$. 
The integer 
\[\cydim(\Delta) := d + 1 - (2 r)\]
we  call the {\em Calabi-Yau dimension} (CY-dimension) of $\Delta$. 
}
\end{definition}

\begin{remark}\ 
{\rm 
\begin{enumerate}
\item The $\Est$-function is symmetric 
\[ \Est(\Delta;u,v) =\Est(\Delta;v,u) \]
and it  satisfies the Poincar{\'e} duality 
\[  (uv)^{\cydimgross(\Delta)} \Est(\Delta;u^{-1},v^{-1} ) =   
\Est(\Delta;u,v); \] 
\item The $\Est$-function  satisfies the following reciprocity-law:
\[\Est(\Delta;u,v)=(-u)^{\cydimgross(\Delta)}\; \Est(\Delta^*;u^{-1},v).\]
Note, that we have by definition $\cydim(\Delta) = \cydim(\Delta^*)$.
\item Let $\deg(\Delta)$ be the degree of the $h^*$-polynomial of 
$\Delta$. Then we have $\cydim(\Delta) = \deg(\Delta) - r$, cf. \cite{Bat07}.
\item The CY-dimension can be arbitrarily negative: If 
$\Delta$ is a lattice pyramid over $\Delta'$, then $\cydim(\Delta) = 
\cydim(\Delta')-1$, 
since the index of $\Delta$ is $1$ plus the index of $\Delta'$.
\end{enumerate}
\label{erez}}
\end{remark}

\smallskip

Here is our conjecture on the $\Est$-function:
\begin{conjecture}{\rm
Let $\Delta$ be a Gorenstein polytope of CY-dimension $n$. 
Then 
$$E := \Est(\Delta;u,v) = \sum_{p,q} (-1)^{p+q} 
h^{p,q}_{\Delta}u^pv^q $$
is a polynomial in $u,v$ of 
degree $2n$ with nonnegative 
integral coefficients $h^{p,q}_{\Delta}$  
(in particular, one has $E = 0$ if 
$n <0$, and  $E$ is constant if  $n=0$).  Moreover, for any $n \geq 1$ 
the polynomial $E$ satisfies the following two conditions
\[ \Est(\Delta;u,0) = (-u)^n  \Est(\Delta;u^{-1},0) \]
and 
\[ \frac{d^2}{du^2}  \Est(\Delta;u,1)|_{u =1} = \frac{n(3n-5)}{12}  
\Est(\Delta;1,1). \]  
These two conditions together with (1) and (2) from \ref{erez}  imply:  
\begin{enumerate} 
\item 
If $n =1$, then for some integer $k$ one has  
\[ \Est(\Delta;u,v)= k (1 - u)(1-v).\] 
\item 
If $n =2$, then  for some integers $k,l$ one has 
\[  \Est(\Delta;u,v) = k (1 + u^2)(1 + v^2) - 2l(u+v)(1+uv) + (20k -16l)uv. \]
\end{enumerate} 
In particular, the number 
$\Est(\Delta;1,1) = 24(k-l)$ is always divisible by $24$.  
\label{e-conj}
}
\end{conjecture}

\begin{remark}{\rm This conjecture is motivated by the expectation 
that the function $\Est(\Delta;u,v)$ ``looks'' 
as a generating function for stringy Hodge 
numbers of a disjoint union of manifolds of dimension $\cydim(\Delta)$
with trivial canonical 
class. The second condition was proved for stringy Hodge numbers 
in a  more  general form in \cite{Bat00}. 
In the case $\cydim(\Delta) =0$ we expect  a finite union of 
points. In the case  $\cydim(\Delta) =1$, it  should be a union of 
elliptic curves.   In the case  $\cydim(\Delta) =2$, one expects 
Hodge numbers of abelian surfaces and $K3$-surfaces.  

In all examples that the authors know the leading coefficient $k$ is either 
zero, or a power of $2$ (see also \ref{alggeonefpart}).
}
\end{remark}

\begin{example}{\rm Let $S_d$ be a standard unimodular simplex of dimension 
$d : = 2r-1$ $(r \in \N$). 
We calculate the $\Est$-function of  
the Gorenstein simplex $\Delta := 2 S_d$ which has index $r$. We have  
$\cydim(\Delta)=0$. 
By Lemma \ref{boundarylemma} 
 the dual Gorenstein polytope $\Delta^*$ is a Gorenstein simplex with only 
unimodular facets. Hence, 
Lemma \ref{pyra} yields:
\[\Est(\Delta;u,v) = (uv)^{-r} \left(\St(\Delta^*,u v) + 
(-u)^{2r} \; \St(\Delta, u^{-1} v)\right).\]
We claim:
\begin{equation}\St(\Delta,t) = t^r = \St(\Delta^*,t).
\tag{$*$}
\label{seq}
\end{equation}
Combining these two equations we get:
\[\Est(\Delta;u,v) = 2.\]

We show (\ref{seq}), by proving that $\St(\Delta,t) \in \{0,t^r\}$ 
for any Gorenstein simplex $\Delta$ of index $r$ and with $\cydim(\Delta)=0$. 
By Prop. \ref{simplexprop} and its proof, to determine $\St(\Delta,t)$ we 
have to compute the sum $\sum_x t^{x_{d+1}}$ over  the lattice points $x$ 
in the interior of the parallelepiped spanned by vertices of $\Delta$. On the 
other hand, in order to determine $h^*_\Delta(t)$ one computes 
the sum   $\sum_x t^{x_{d+1}}$ over {\em all} 
lattice points $x$  in the paralellepiped. 
By Remark \ref{erez}(3), the degree of $h^*_\Delta$ equals the index 
$r$. 
Hence, there is only  one interior lattice point $x$ 
in the parallelepiped such that  $x_{d+1} = r$.

In this case, $x$ is the unique interior lattice point 
in $r \Delta = (2 r) S_d$ which 
is given as $\sum_{v \in \V(\Delta)} \frac{1}{2} v$. Hence $\St(\Delta,t)=t^r$. 
In the same way, this holds 
for $\Delta^*$. This proves (\ref{seq}).
}
\end{example}

\smallskip

From Remark \ref{erez}(4) we see that using lattice pyramids we may 
construct Gorenstein polytopes of  negative CY-dimension. 
In this case, the $\Est$-function should vanish. However, this vanishing 
 holds for any lattice pyramid (even of positive CY-dimension):

\begin{proposition}
Let $\Delta$ be a Gorenstein polytope. If $\Delta$ is a lattice 
pyramid, then $\Est(\Delta;u,v) = 0$.
\label{pyramid}
\end{proposition}

\begin{proof}

Let $\Delta$ be a lattice pyramid over a facet $\Delta'$ with apex $v$. Then 
also $\Delta^*$ is a lattice pyramid over the facet $v^*$ with apex 
$(\Delta')^*$, 
for instance, use Prop. \ref{char}. 
Thus, for any face $F$ of $\Delta$, $F$ or $F^*$ is a lattice pyramid. 
Hence, $\St(F,u^{-1} v) \, \St(F^*, u v) = 0$ by Lemma \ref{pyra}. From this 
the statement follows.
\end{proof}

\smallskip

\subsection{The $B$-polynomial}

Originally, in \cite{BB96b} so-called {\em $B$-polynomials} were used 
to calculate 
Hodge numbers, see also \cite[Remark 5.8]{KRS03} and \cite{Haa00}. 
In \cite[Lemma 11.4]{BM03} Borisov and Mavlyutov gave the following 
definition in terms of $g$-polynomials:

\begin{definition}{\rm 
Let $\calP$ be an Eulerian poset of rank $d$, with minimal element 
$\hat{0}$ and maximal element $\hat{1}$. 
Then \[B(\calP;u,v) := \sum_{\hat{0} \leq x \leq \hat{1}} \, 
(-u)^{\rk(\hat{1}) -\rk(x)}\; 
g_{[x,\hat 1]^*}(u^{-1} v)\;  g_{[\hat{0},x]}(uv).\]
}
\end{definition}

To see, how the $B$-polynomial naturally appears, 
we separately collect all $g$-polynomials in the definition of 
the $\Est$-function 
into one polynomial in two variables: Let $\Delta$ be a 
Gorenstein polytope of index $r$. 
Then $\Est(\Delta;u,v)$ equals 
\[(uv)^{-r} \sum_{F, F'}  
(-1)^{\dim((F')^*)+1} u^{\dim(F)+1} \; h^*_F(u^{-1} v) \; h^*_{F'}(u v) \; 
B([F,(F')^*];u,v),\]
where the sum is over all pairs of faces $\emptyset \leq F \leq \Delta$, 
$\emptyset \leq F' \leq F^*$.

\smallskip

There is an important identity on $g$-polynomials, called Stanley's 
convolution, cf. \cite[Cor.8.3]{Sta92}, 
that can be stated in the following form, see \cite[Lemma 11.2]{BM03},

\begin{lemma}
$B([F,(F')^*];1,v)$ equals $0$, if $F \not= (F')^*$, and $1$, 
otherwise.
\end{lemma}

As a corollary we get
\[\Est(\Delta;1,v) = v^{-r} \sum_{\emptyset \leq F \leq \Delta} 
(-1)^{\dim(F)+1} \;h^*_F(v)\; h^*_{F^*}(v),\]
\begin{equation}
\Est(\Delta;1,1) = \sum_{\emptyset \leq F \leq \Delta} (-1)^{\dim(F)+1} 
\;\Vol(F)\; \Vol(F^*),
\tag{$**$}
\label{comby}
\end{equation}
where $\Vol$ is the lattice volume (with convention $\Vol(\emptyset) := 1$). 

\begin{remark}{\rm 
When $\Delta$ is a three-dimensional reflexive polytope, Theorem \ref{hodge} below yields that 
$\Est(\Delta;1,1)$ equals the topological Euler-characteristic of a $K3$ surface. Hence, 
$\Est(\Delta;1,1) = 24$. However, there is still no direct combinatorial proof for 
this equation, cf. \cite{Haa05}.
}
\end{remark}

\smallskip

\subsection{Stringy Hodge numbers} 

The main motivation for the definition of the $\Est$-function 
and the Calabi-Yau dimension, and 
for the related conjectures stems from the following theorem:

\begin{theorem}\cite[Thm.7.2]{BM03}
Let $\sigma$ be a $(d+r)$-dimensional reflexive Gorenstein cone of index $r$ 
with support $\Dt \subseteq \MRt$, where $\Dt$ is a Cayley polytope 
of $r$ lattice polytopes. 
Let $Y$ be a generic complete intersection Calabi-Yau defined by these $r$ 
equations, hence $\dim(Y) = d-r = \cydim(\Dt)$. 
Then 
\[\Est(\Dt;u,v) = \Est(Y;u,v) = \sum_{p,q}(-1)^{p+q} \; 
\hst^{p,q}(Y) \; u^p \; v^q\]
is a polynomial in $u,v$, called {\em $\Est$-polynomial} 
of $Y$, where the coefficients 
$\hst^{p,q}$ are the so-called {\em stringy Hodge numbers} 
of $Y$, in the sense of \cite{BD96,Bat98}.
\label{hodge}
\end{theorem}

This is a simplified version of the main result in \cite{BB96b}. 
Note, that the combinatorial identity (\ref{comby}) 
gives the topological Euler-characteristic of the 
Calabi-Yau $Y$. From this theorem 
we can derive the following result using algebraic geometry 
(for $n=2$, $Y$ is a $K3$-surface) 
and the reciprocity-law (Remark \ref{erez}(2)):

\begin{corollary}
Conjecture \ref{e-conj} holds, if at least one of the two Gorenstein 
polytopes $\Delta$ or $\Delta^*$ of index $r$ is a Cayley polytope of 
length $r$.
\end{corollary}

\begin{example} 
{\rm Here is an example for a Gorenstein polytope $\Delta$ that does not 
satisfy the 
assumption of the previous corollary. Let $\Delta_1 := 2 S_3$, where $S_3$ 
is the unimodular $3$-simplex, and let $\sigma_1$ be the 
reflexive Gorenstein cone of index two with support $\Delta_1$. 
We also define $\Delta_2 := \Delta_1^*$, and $\sigma_2 := \sigma_1^\dual$.
 Then by Lemma \ref{directsum} below, the support $\Delta$ of the 
direct sum $\sigma_1 \oplus \sigma_2$ is a $7$-dimensional Gorenstein 
polytope of index $4$. The dual Gorenstein polytope $\Delta^*$ is the support 
of $\sigma_2 \oplus \sigma_1$. Since $\Delta_2$ does not contain a 
special $1$-simplex, both $\Delta$ and $\Delta^*$ do not contain special 
$3$-simplices, hence by Prop. \ref{ath} 
both Gorenstein polytopes of index $4$ are not Cayley polytopes of length 
$4$.} 
\end{example} 

The proof of the following lemma is left to the reader.

\begin{lemma}
Let $\sigma_1 \subseteq \Mq_\R$, $\sigma_2 \subseteq \Mq'_\R$ be two reflexive Gorenstein cones of indices 
$r_1$, respectively $r_2$. Then the direct sum 
$\sigma_1 \oplus \sigma_2 \subseteq (\Mq \oplus_\Z \Mq')_\R$ is a reflexive Gorenstein cone of index 
$r_1 + r_2$, with dual cone $\sigma_1^\dual \oplus \sigma_2^\dual$. 
Moreover, the support of $\sigma_1 \oplus \sigma_2$ contains a special 
$(r_1+r_2-1)$-simplex if and only if 
the support of $\sigma_1$ contains a special $(r_1-1)$-simplex and the 
support of $\sigma_2$ contains a special $(r_2-1)$-simplex.
\label{directsum}
\end{lemma}

If in Theorem \ref{hodge} the cone $\sigma$ is associated to a nef-partition,
 then we get from Remark \ref{erez}(2) 
for the $\Est$-polynomial of $Y$
\[\Est(Y;u,v)=(-u)^{d-r} \Est(Y^*;u^{-1},v),\]
where $Y^*$ is a $(d-r)$-dimensional complete intersection 
Calabi-Yau defined by the dual nef-partition. 
In particular, stringy Hodge numbers satisfy the following mirror 
symmetry property:
\[\hst^{p,q}(Y) = \hst^{d-r-p,q}(Y^*) \quad \text{ for } 0 \leq p,q \leq d-r.\]

\smallskip

\subsection{Boundedness of Calabi-Yau $n$-folds}

It is still an open question  whether there is only a 
finite number of topologically different Calabi-Yau $n$-folds ($n \geq 3$). 
We can look at this question from a combinatorial point of view and ask 
another one:

\begin{question}{\rm Let $n$ be fixed. 
Is there up to scalar only a finite number of polynomials in $u,v$ that occur 
as $\Est$-polynomials of Gorenstein polytopes $\Delta$ 
such that  $\cydim(\Delta)$ equals a fixed number $n$?
\label{conj}
}
\end{question}

\begin{remark}{\rm 
In particular, a positive answer would imply that 
the set of stringy Hodge numbers associated to irreducible 
Calabi-Yau complete intersections $Y$ defined by a nef-partition is bounded. 
This implication follows from the fact that by Poincar{\'e} duality 
the leading coefficient $k$ of $\Est(Y;u,v)$ equals the 
constant coefficient which is $1$, 
because $k = h^{0,0}(Y) = h^0(O_Y)$ 
is  the number of irreducible components of $Y$. 
We observe that in general the number $k$ can be arbitrary 
large powers of $2$.  Indeed, let $\Delta_1, 
\ldots, \Delta_r$ be a nef-partition of length $r$ and 
let $E$ be  $\Est$-polynomial of the corresponding Calabi-Yau 
complete intersection $Y$ of dimension $d-r$. 
Then $r+1$ polytopes 
$$(\Delta_1 \times 0), \ldots, (\Delta_r \times 0), (0 
\times [-1,1]) \subseteq \MR \oplus \R$$ 
define a  reducible nef-partition of length $r+1$ of the same  
CY-dimension $d-r$. 
The corresponding Calabi-Yau  complete intersection $Y'$ is a product 
of $Y$ with  two points.   
Therefore $\Est$-polynomial of $Y'$ equals $2 E$. By repeating the same  
procedure, we get  $\Est$-polynomials with arbitrary large 
leading coefficient $k$.
\label{alggeonefpart}
}
\end{remark}

\begin{remark} 
{\rm There are various finiteness results on toric Fano varieties 
\cite{Bat82,BorBor93}. However these results are not sufficient 
to prove  finiteness of $\Est$-polynomials of Gorenstein 
polytopes of fixed Calabi-Yau dimension. We are still optimistic 
and want to illustrate our hope by  an example of $d$-dimensional 
Gorenstein simplices $S(\omega)$  
of index $r$ defined as 
\[ S(\omega) := \{ (x_0,x_1 \ldots, x_d) \in \R_{\geq0}^{d+1} \; : \; 
w_0 x_0 + w_1x_1 + \cdots + w_d x_d = w \}, \]
where $w, w_0, w_1, \ldots, w_d$ are positive integers having the property
$w_i |w$ $\forall i$. Let $k_i: = w/w_i$ $( 0 \leq i \leq d)$. The fact 
that  $S(\omega)$ is a Gorenstein simplex of index $r$ is equivalent 
to the equation
\[ \frac{1}{k_0} +  \frac{1}{k_1} + \cdots +  \frac{1}{k_d} = r. \]
If $k_i =1$ for some $i$, then   $S(\omega)$ is a lattice pyramid over
a $(d-1)$-dimensional simplex and, by \ref{pyramid}, the $\Est$-polynomial
of  $S(\omega)$ is $0$. Therefore it is enough to consider only the case 
when $k_i \geq 2$ $( 0 \leq i \leq d)$. Let $s$ be the number of indices 
$i$ such that $k_i \geq 3$. Then we get
from the above equation  
\[  r = \frac{1}{k_0} +  \frac{1}{k_1} + \cdots +  \frac{1}{k_d}  \leq 
\frac{ d+1 -s}{2} + \frac{s}{3}, \] 
or $3(d +1 -2r) \geq s$. By definition  
\[  \cydim(S(\omega)) = d+1 -2r. \]
So we obtain $s \leq 3 \cydim(S(\omega))$. This implies that if  $S(\omega)$
has  a fixed Calabi-Yau dimension $n = d+1 -2r$, then 
there exist only 
finitely many possibilities for integers $k_i$ which are different 
from $2$. In order to compute  $\Est(S(\omega);u,v)$ we observe 
that the simplex $S(\omega)$ determines a Landau-Ginzburg model 
with LG-potential of  Fermat type: 
\[ F(z) := z_0^{k_0} +  z_1^{k_1} + \cdots +   z_d^{k_d}. \]
The corresponding formula for the $\Est$-polynomial is expected to coincide
with the well-known formula of Vafa \cite{Va89}.  The latter one 
does not change if one adds any quadratic terms to $F(z)$. Therefore we 
come to finitely many  $\Est$-polynomials for a fixed value 
of  $\cydim(S(\omega))$. 
} 
\end{remark}

\section{Nef-partitions and  Homological Mirror Symmetry}

Let $P$ and $P^*$ be dual to each other Gorenstein polytopes 
of index $r>1$. Example \ref{cube} shows that $P$ and $P^*$ may contain 
many different special $(r-1)$-simplices. 
Consider the case when a Gorenstein polytope 
$P$ contains two different special $(r-1)$-simplices $S'$ and $S''$. 
Assume that the dual Gorenstein polytope $P^*$ also contains  a special 
$(s-1)$-simplex $S$.  
Then the special $(r-1)$-simplex $S$ defines a Cayley structure on 
$P$, i.e., some lattice polytopes $\Delta_1, \ldots, \Delta_r$ 
such that $P \cong \Delta_1 * \cdots * \Delta_r$. Now two different 
special $(r-1)$-simplices $S', S'' \subset P$ define two collections 
of lattice points $p_i', p_i'' \in \Delta_i$ ($ 1\leq i \leq r)$ such that 
\[ m = p_1' + \cdots + p_r' = p_1'' + \cdots + p_r'' \]
is the unique interior lattice point in the reflexive polytope
$\Delta_1 + \cdots + \Delta_r$. We remark that two Cayley structures of 
the dual Gorenstein polytope $P^*$ corresponding to $S'$ and $S''$ 
may define very different polytopes $\nabla_1, \ldots, \nabla_r$.

\begin{example}{\rm
We consider below a picture of two centered nef-partitions corresponding to 
polytopes $\Delta_1', \Delta_2'$ and $\Delta_1'', \Delta_2''$ which are 
the same up to a affine translation, i.e.    
$\Delta_1'' = \Delta_1' - (0,1)$ and 
$\Delta_2'' = \Delta_2' - (0,-1)$:\medskip

{\centerline{
\psfrag{a}{$\Delta_1'$}
\psfrag{b}{$\Delta_2'$}
\psfrag{c}{$\nabla_2'$}
\psfrag{d}{$\nabla_1'$}
\psfrag{e}{$\Delta_2''$}
\psfrag{f}{$\Delta_1''$}
\psfrag{g}{$\nabla_1''$}
\psfrag{h}{$\nabla_2''$}
\includegraphics{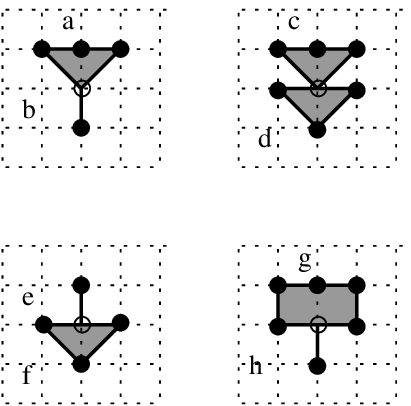}}}
  \label{fig2}
We see that  the corresponding polytopes  $\nabla_1', \nabla_2'$ and 
$\nabla_1'', \nabla_2''$ are not isomorphic via affine translations by 
lattice vectors. } 
\end{example}

It is easy to see  that, if two nef-partitions consist of 
collections of lattice polytopes 
$\{ \Delta_1', \ldots, \Delta_r'\}$ and 
$\{ \Delta_1'', \ldots, \Delta_r''\}$ such that 
$\Delta_i'' = \Delta_i' + q_i$  for $q_i \in M$ $( 1\leq i \leq r)$, 
then the corresponding generic Calabi-Yau complete intersections 
$\overline{X}_{\Delta_i'} $ and $\overline{X}_{\Delta_i''}$ are isomorphic. 
However, there is no natural isomorphism between the generic mirror 
 Calabi-Yau complete intersections 
$\overline{X}_{\nabla_i'} $ and $\overline{X}_{\nabla_i''}$. 
Using the formula for $E_{\rm st}$-polynomial of $P^*$, we can only say that 
the stringy Hodge numbers of 
$\overline{X}_{\nabla_i'} $ and $\overline{X}_{\nabla_i''}$ are the same.

\begin{question} 
Are the mirror Calabi-Yau complete intersections 
$\overline{X}_{\nabla_i'} $ and $\overline{X}_{\nabla_i''}$
birationally isomorphic?
\end{question} 

Even if in general  the answer to this question will be negative, we can  
use the approach of the Homological Mirror Symmetry and formulate 
the following conjecture: 

\begin{conjecture} 
There exists a equivalence (of Fourier-Mukai type) between the derived 
category of coherent sheaves on 
two Calabi-Yau complete intersections 
$\overline{X}_{\nabla_i'} $ and $\overline{X}_{\nabla_i''}$. 
\end{conjecture} 

\begin{remark} 
{\rm This conjecture needs some comments, because Calabi-Yau varieties 
$\overline{X}_{\nabla_i'} $ and $\overline{X}_{\nabla_i''}$ are usually 
singular. However, using toric methods, one can always construct 
their partial crepant desingularizations 
$\widehat{X}_{\nabla_i'} $ and $\widehat{X}_{\nabla_i''}$ which are 
Calabi-Yau varieties having at worst terminal abelian 
quotient singularities. As in the toroidal case considered by Kawamata 
in \cite{Ka05} the equivalence class of stacky derived categories of 
$\widehat{X}_{\nabla_i'} $ and $\widehat{X}_{\nabla_i''}$ should be 
independent of the choice of the  partial crepant desingularizations. 
So we expect the existence of a  Fourier-Mukai transform which establishes 
an equivalence between the  stacky derived categories of 
$\widehat{X}_{\nabla_i'} $ and $\widehat{X}_{\nabla_i''}$. 
} 
\end{remark} 

\section{Operations with nef-partitions}

In this concluding section we give as an addendum to Section 3 
a concise survey on how to modify, project and decompose nef-partitions. 

\smallskip

Throughout, let $\Delta_1, \ldots, \Delta_r \subseteq \MR$ and 
$\nabla_1, \ldots, \nabla_r \subseteq \NR$ 
be two dual to each other centered proper nef-partitions.

\smallskip \subsection{Collecting} 
By collecting lattice polytopes in a nef-partition we obtain 
another nef-partition. 
For this, we simply replace some lattice polytopes in a 
nef-partition by their Minkowski sum. This yields a shorter 
nef-partition, which is by Cor. \ref{nabcor} dual to the 
nef-partition we get by replacing the respective lattice polytopes in 
the dual nef-partition by their convex hull:

\begin{definition}\ { \rm 
For $I \subseteq \{1, \ldots, r\}$ we define $\nabla_I := 
\conv(\nabla_i \,:\, i \in I)$ and $\nabla^I := \sum_{i \in I} \nabla_i$. 
In the same way we define $\Delta_I$ and $\Delta^I$.
}
\end{definition}

\begin{proposition}
Let $I_1, \ldots, I_l$ be a partition of $\{1, \ldots, r\}$ 
into non-empty subsets. Then 
$\Delta^{I_1}, \cdots, \Delta^{I_l}$ and  $\nabla_{I_1}, \ldots, 
\nabla_{I_j}$ are dual nef-partitions. 
\label{regroup}
\end{proposition}

From Prop. \ref{nefdualdef} we get:

\begin{corollary}
For $I \subseteq \{1, \ldots, r\}$ we have
\[\Delta^I = \{x \in \Delta_1 + \cdots + \Delta_r \,:\, 
\pro{x}{\nabla_j} \geq 0 \;\forall\, j \not\in I\},\]
\[\nabla_I = \{y \in \conv(\nabla_1, \ldots, \nabla_r) \,:\, 
\pro{\Delta_j}{y} \geq 0 \;\forall\, j \not\in I\}.\]
\label{suppfct}
\end{corollary}

\begin{remark}{\rm Let us explain how Prop. \ref{regroup} 
relates to the results of Section~2: 
More generally, let $\Delta_1, \ldots, \Delta_r$ be arbitrary 
lattice polytopes with $\Delta_1 + \cdots + \Delta_r$ reflexive. 
Then by Theorem \ref{sum}, 
$P := \Delta_1 * \cdots * \Delta_r$ is a Gorenstein Cayley polytope 
of index $r$. The dual Gorenstein polytope $P^*$ contains a special 
$(r-1)$-simplex 
with vertices $v_1, \ldots, v_r$ (see Prop. \ref{ath}). 
Theorem \ref{brtheo} shows that projecting $P^*$ along the affine 
hull of this special simplex yields the reflexive polytope 
$(\Delta_1 + \cdots + \Delta_r)^*$. 

Now, let $I_1, \ldots, I_l$ be a partition of $\{1, \ldots, r\}$ 
into non-empty subsets. 
We may assume $I_1 = \{1, \ldots, r'\} \subseteq \{1, \ldots, r\}$. 
In the same way as in the proof 
of Theorem \ref{sum} and Theorem \ref{brtheo}, projecting $P^*$ along 
$\aff(v_1, \ldots, v_{r'})$ yields a Gorenstein polytope dual to 
the Gorenstein polytope $\Delta^{I_1} * \Delta_{r'+1} * \cdots * 
\Delta_r$ of index $r+1-r'$. 
After $l$ such projection steps, we finally get a Gorenstein polytope
 of index $l$ dual to $\Delta^{I_1} * \Delta^{I_2} * \cdots * \Delta^{I_l}$.

In particular, when $\Delta_1, \ldots, \Delta_r$ and $\nabla_1, \ldots, 
\nabla_r$ are dual to each other centered nef-partitions, we see that 
$P^* = \nabla_1 * \cdots * \nabla_r$, and the projection along $\aff(v_1, 
\ldots, v_{r'})$ yields 
$\nabla_{I_1} * \nabla_{r'+1} * \cdots * \nabla_r$ as the Gorenstein 
polytope dual to 
$\Delta^{I_1} * \Delta_{r'+1} * \cdots * \Delta_r$. Hence, eventually we get 
$\nabla_{I_1} * \nabla_{I_2} * \cdots * \nabla_{I_l}$ as the Gorenstein 
polytope of index $l$ dual 
 $\Delta^{I_1} * \Delta^{I_2} * \cdots * \Delta^{I_l}$.
}
\end{remark}

\smallskip \subsection{Projecting}

Here, we show that projecting along a lower-dimensional lattice polytope 
in a nef-partition still yields a nef-partition. 

\smallskip

The following definition is convenient:

\begin{definition}{\rm 
Let $A \subseteq \MR$ be a set. Then 
\[A^\perp := \lin(A)^\perp \subseteq \NR.\]
}\end{definition}

We start with a general observation:

\begin{lemma}
Let $P,Q$ be lattice polytopes with $\dim(P) < \dim(P+Q)$. We denote by 
$\pi$ the projection along $\aff(P)$. 
Then preimages of facets of $\pi(P+Q)$ are facets of $P+Q$.
\label{minky}
\end{lemma}

\begin{proof}

By translation we may assume $\lin(P) = \aff(P)$. 
Let $G$ be a facet of $\pi(P+Q)$. Then the preimage $F$ of $G$ is a 
Minkowski sum $F = P + F_Q$ for a face $F_Q$ of $Q$. 
Since $\pi(P+Q) = \pi(Q)$, we easily see that $\dim(F) = \dim(P) + 
\dim(G)$. Now, the statement follows 
from $\dim(G) = \dim(P+Q) - \dim(P) - 1$.
\end{proof}

\begin{corollary}
Let $P,Q$ be lattice polytopes with $P+Q$ reflexive. Then the projection 
$\pi(P+Q)$ along $\aff(P)$ is a reflexive polytope.
\label{minkyrefl}
\end{corollary}

\begin{proof}

We may suppose $\dim(P) < \dim(P+Q)$. 
Let $m$ be the unique interior lattice point of $P+Q$. By considering 
$P-m,Q$ instead of $P,Q$ 
we may assume that $P+Q$ is reflexive with respect to $0$.
Let $p$ be any lattice point of $P$. Then by considering $P-p,Q+p$ 
instead of $P,Q$ we may also assume that $0 \in P$. 
Hence $\aff(P) = \lin(P)$. So $\pi(P+Q)$ has $0$ as an interior lattice 
point. Now, Lemma \ref{minky} yields that 
any vertex of $(\pi(P+Q))^* = (P+Q)^* \cap P^\perp$ is also a vertex of 
$(P+Q)^*$. Since $(P+Q)^*$ is a lattice polytope, 
$(\pi(P+Q))^*$ is also a lattice polytope.
\end{proof}

Now, we deal with our special situation of dual to each other centered 
nef-partitions $\Delta_1, \ldots, \Delta_r$ and $\nabla_1, \ldots, \nabla_r$:

\begin{proposition}
Let $\{1, \ldots, r\} = I \sqcup J$ with $I \not= \emptyset \not= J$, 
such that $\lin(\Delta^J)$ has dimension strictly lower than $\MR$. We 
denote by $\pi$ the projection map along $\lin(\Delta^J)$. 
Then $\pi(\Delta_1), \ldots, \pi(\Delta_r)$ is a centered nef-partition.

The dual reflexive polytope (in $(\Delta^J)^\perp$)
\[(\pi(\Delta_1) + \cdots +  \pi(\Delta_r))^* = \nabla_I \cap (\Delta^J)^\perp\]
is the smallest face $F$ of $\nabla_I$ containing $0$.

The dual nef-partition is given as $\nabla_1 \cap F, \ldots, \nabla_r \cap F$.
\label{projtheo}
\end{proposition}

\begin{proof}

By Cor. \ref{suppfct}, we have 
\[(\pos(\Delta^J))^\dual = \pos(\nabla_I),\]
and 

\begin{center}
$\begin{array}{ll}
&(\pi(\Delta_1) + \cdots +  \pi(\Delta_r))^* = (\Delta_1 + \cdots + 
\Delta_r)^* \cap (\Delta^J)^\perp =\\
= &\conv(\nabla_I, \nabla_J) \cap (\Delta^J)^\perp = \nabla_I \cap 
(\Delta^J)^\perp = F.
\end{array}$
\end{center}

The other statements are clear.
\end{proof}

\begin{corollary}
Let $\{1, \ldots, r\} = I \sqcup J$ with $I \not= \emptyset \not= J$. 
Let $F$ be the smallest face $F$ of $\nabla_I$ containing $0$. 
Then
\[\dim(F) + \dim(\Delta^J) = \dim(\MR).\]
\end{corollary}

In particular, $\dim(\Delta_1) = \dim(\MR)$ iff $0$ is a vertex of 
$\conv(\nabla_2, \ldots, \nabla_r)$.

\smallskip \subsection{Decomposing}

In \cite{BB96a} Borisov and the first author introduced the notion of an 
irreducible nef-partition:

\begin{definition}{\rm 
Let $\emptyset \not= I \subseteq \{1, \ldots, r\}$ such that $\Delta^I$ is
 a reflexive polytope (in its linear span). By Prop. \ref{projtheo} it is enough to assume 
that $\Delta^I$ contains $0$ in its relative interior. 
We say $I$ is {\em irreducible}, if $I$ is minimal with this property. 

The proper nef-partition $\Delta_1, \ldots, \Delta_r$ is called 
{\em irreducible}, if $\{1, \ldots, r\}$ is irreducible.
}
\end{definition}

Theorem~5.8 of \cite{BB96a} states that any (proper) nef-partition has 
a unique {\em decomposition} into 
irreducible (proper) nef-partitions, i.e., there are unique irreducible 
subsets $I_1, \ldots, I_l$ that form a partition of 
$\{1, \ldots, r\}$ such that 
\[\Delta_1 + \cdots + \Delta_r = \Delta^{I_1} + \cdots + 
\Delta^{I_l},\quad \lin(\Delta^{I_1}) \oplus \cdots \oplus \lin(\Delta^{I_l}) = 
\MR.\]
This is called {\em semi-simplicity of nef-partitions}. 
Here, we recall the proof of the existence of such a decomposition and 
give a generalization in Prop.~\ref{cancel}. We 
also prove in Prop.~\ref{unique} the uniqueness of this decomposition,
 which was missing in \cite{BB96a}. 
Finally, we determine the maximal length of a proper nef-partition in 
Prop.~\ref{length}.

\smallskip

The existence of a decomposition into irreducible nef-partitions is 
an immediate consequence of the following result, 
which was shown in \cite[Thm.5.8]{BB96a} (we include a proof for the 
convenience of the reader):

\begin{proposition}
Let $\{1, \ldots, r\} = I \sqcup J$ with $I \not= \emptyset \not= J$. 
If $\Delta^I$ is reflexive, then 
also $\Delta^J$ is reflexive, and $\lin(\Delta^I) \oplus_\R \lin(\Delta^J) = 
\MR$.
\label{exprop}
\end{proposition}

\begin{proof}

By Prop.~\ref{regroup} $\Delta^I,\Delta^J$ and $\nabla_I, \nabla_J$ are 
dual nef-partitions. 
Since $\Delta^I$ contains $0$ in its relative interior, Prop. 
\ref{projtheo} 
yields $\Delta^I \subseteq U := \nabla_J^\perp$. Moreover, 
$\Delta^J \subseteq \sigma := (\pos(\nabla_I))^\dual$ because of 
$\pro{\Delta^J}{\nabla_I} \geq 0$. Furthermore, $U \cap \sigma = \{0\}$, since 
$\pro{x}{\nabla_I+\nabla_J} \geq 0$ is only possibly for $x=0$. 
We show that $\sigma$ is actually a linear subspace of $\MR$. 
This already yields $U \oplus \sigma = \lin(\Delta^I) \oplus \lin(\Delta^J) 
= \MR$, 
and by projecting along $U = \lin(\Delta^I)$ we get that $\Delta^J$ is reflexive by Prop. \ref{projtheo}.

Let $x \in \sigma$. We have to show $-x \in \sigma$. 
Let $\rho$ be the projection along~$U$. As $0$ is in the interior of 
$\Delta^I+\Delta^J$, and $\rho(\Delta^I) = \{0\}$, 
we get $\rho(\sigma) \supseteq \rho(\pos(\Delta^J)) = \rho(\MR)$. 
Therefore there is some $x' \in \sigma$ with 
$\rho(x') = \rho(-x)$. This implies $x'+x \in U \cap \sigma = \{0\}$, 
thus $-x = x' \in \sigma$.
\end{proof}

As was observed in \cite[Example 5.5]{BB96a}, note that the direct sum 
in Prop. \ref{exprop} 
is in general only a splitting of vector spaces, {\em not} of lattices.

\begin{corollary} 
Let $\{1, \ldots, r\} = I \sqcup J$ with $I \not= \emptyset \not= J$. 
If one of the lattice polytopes in 
$\{\Delta^I, \Delta_I, \Delta^J, \Delta_J, \nabla^I, \nabla_I, \nabla^J, 
\nabla_J\}$ contains $0$ in its relative interior, then 
anyone does, and, moreover, each one is a reflexive polytope. 

In this case: $(\Delta_i)_{i \in I}$ is a nef-partition with dual 
nef-partition $(\nabla_i)_{i \in I}$, 
with respect to the dual vector spaces $\nabla_J^\perp \subseteq \MR$ and 
$\Delta_J^\perp \subseteq \NR$.
\label{corex}
\end{corollary}

\begin{proof}

By Prop.~\ref{regroup}, $\Delta^I, \Delta^J$ and $\nabla_I, \nabla_J$ are
 dual nef-partitions. The same is true for 
$\nabla^I, \nabla^J$ and $\Delta_I, \Delta_J$. Hence, Prop. \ref{projtheo}, 
implies that if some of these $8$ polytopes contains $0$ in its relative 
interior, 
then it is necessarily a reflexive polytope. Since $\Delta^I$ has $0$ in 
its relative interior if and only if $\Delta_I$ has $0$ in its relative 
interior, 
we may assume that $\Delta^I$ is reflexive. Prop. \ref{exprop} yields that 
also $\Delta^J$ is reflexive. Now, Cor. \ref{suppfct} implies
\[\nabla_I = \{y \in (\pos(\Delta^J))^\dual \,:\, \pro{\Delta^I}{y} \geq -1\} = 
\{y \in (\Delta_J)^\perp \,:\, \pro{\Delta^I}{y} \geq -1\},\]
since $\Delta_J$ has $0$ in its relative interior. From this, we get 
$\nabla_I = (\Delta^I)^*$. Hence, $\nabla_I$ is reflexive. In the same
 way we see 
$\nabla_J = (\Delta^J)^*$, so also $\nabla_J$ is reflexive. 
Now, the last statement follows merely from the definition of dual 
nef-partitions.
\end{proof}

There is the following generalization of Prop. \ref{exprop} 
(this result is joint work with Christian Haase):

\begin{proposition} 
Let $P,Q$ be lattice polytopes with $P+Q$ reflexive. If $P$ contains a 
lattice point in its relative interior, 
then $P$ and $Q$ are reflexive polytopes, and $P,Q$ form a nef-partition.
\label{cancel}
\end{proposition}

\begin{proof}

As at the beginning of the proof of Cor. \ref{minkyrefl} we may assume that 
$0$ is an interior lattice point of $P$ and $P+Q$.

We define $Q' := Q \cap \lin(P)$. Since $0 \in P+Q$, there exist $x \in P, y 
\in Q$ with $x+y = 0$. Hence $Q' \not= \emptyset$. 
Assume there is some $q \in Q'$, $q \not= 0$. 
Since $0$ is in the relative interior of $P$, there is a linear form $u$ in 
the dual space of $\lin(P)$ such that 
$\pro{u}{q} < 0$, and $F := \{x \in P \,:\, \pro{u}{x} = -1\}$ is a facet of
 $P$. Let $y \in Q'$ 
with $\pro{u}{y} = \min_{x \in Q'} \pro{u}{x} =: c$. We have 
$c \leq \pro{u}{q} < 0$. 
Hence, $G := \{x \in P+Q' \,:\, \pro{u}{x} = -1 + c\}$ is a face of 
$ P + Q' = (P+Q) \cap \lin(P)$. Since $0 \not\in G$, $F + y \subseteq G$ and 
$\dim(F) = \dim(\lin(P)) - 1$, we see that $G$ is a facet of $(P+Q) 
\cap \lin(P)$, parallel to $F$. 
Since $P+Q$ is reflexive, there is an integral linearform $w$ in the 
dual space of $\lin(P)$ such that $\pro{w}{G} = -1$. 
We have $u = \lambda w$ for $\lambda := 1 - c > 1$. Let $v$ be a vertex 
of $F$, hence 
$\pro{w}{v} \in \Z$. This implies $\lambda = \pro{u}{v}/\pro{w}{v}
 = -1/\pro{w}{v} \leq 1$, a contradiction.

We have shown that $Q \cap \lin(P) = \{0\}$. In particular $0 \in Q$, 
thus $P,Q$ is by definition a centered nef-partition. Now, 
apply Cor. \ref{corex} and Prop. \ref{exprop}.
\end{proof}

\begin{example}{\rm Let $P := \conv((1,0,0),(0,1,0),(-1,-1,0))$ and $Q := 
\conv((1,0,-1),(0,1,1))$. 
Then $P+Q$ has only one interior lattice point, $P$ is reflexive, however 
$Q$ has no interior lattice point. 
This example, due to Christian Haase, shows that in Prop. \ref{cancel} 
the assumption that $P+Q$ is reflexive 
cannot be weakened by the condition that $P+Q$ has only one interior 
lattice point.
\label{without}
}
\end{example}

\smallskip

The uniqueness of a decomposition into irreducible nef-partitions is 
a corollary of the following result:

\begin{proposition} Any two irreducible subsets of $\{1, \ldots, r\}$ 
are disjoint.
\label{unique}
\end{proposition}

\begin{proof}
Assume there exist irreducible subsets $I,I'$ of $\{1, \ldots, r\}$ 
with $I \not= I'$ and $A := I \cap I' \not=\emptyset$. 
Thus, we have the following situation: 
$\{1, \ldots, r\} = I \sqcup J = I' \sqcup J'$, $I = A \sqcup B$, 
$J = C \sqcup D$, $I' = A \sqcup C$, $J' = B \sqcup D$.

By Prop. \ref{exprop} we have $\lin(\Delta^I) \cap \lin(\Delta^J) = \{0\}$, thus 
$\lin(\Delta^A) \cap \lin(\Delta^C) = \{0\}$. 
Hence, since $\Delta^{I'} = \Delta^A + \Delta^C$ contains $0$ in its relative interior, 
also $\Delta^A$ contains $0$ in its relative interior. This contradicts $I$ being irreducible.
\end{proof}

\smallskip

Finally, we examine nef-partitions of largest length:

\begin{proposition}
A proper nef-partition has length $r \leq 2 \dim(\MR)$, with equality 
only if the convex hull is 
combinatorially a crosspolytope, i.e., the dual of a cube.
\label{length}
\end{proposition}

\begin{proof}

An irreducible proper nef-partition has length $r \leq \dim(\MR)+1$, 
see \cite[Cor.3.5]{BB96a} or \cite[Lemma 2.6]{HZ05}. 
Now, a proper nef-partition of length $r$ decomposes into, say $l$, 
irreducible nef-partitions. 
This yields 
\[r \leq \dim(\MR) + l \leq 2 \dim(\MR),\]
where complete equality holds if only if any irreducible nef-partition 
has dimension one and length two.
\end{proof}

\end{document}